\documentclass[lettersize,journal]{IEEEtran}
\usepackage{amsmath,amsfonts,amssymb,cases,mathdots,color,colordvi,url,tikz}
\usepackage{hyperref}
\usepackage{bm}
\usepackage{bbm}
\usepackage{pgfpages}
\usepackage{algorithmic}
\usepackage{algorithm}
\usepackage{array}
\usepackage[caption=false,font=small,labelfont=rm,textfont=rm]{subfig}
\usepackage{textcomp}
\usepackage{pgfpages}
\usepackage{stfloats}
\usepackage{url}
\usepackage{verbatim}
\usepackage{graphicx}
\usepackage{cite}
\usepackage{multirow}
\usepackage{caption}
\usepackage{enumitem}
\newtheorem{proposition}{Proposition}
\newtheorem{theorem}[proposition]{Theorem}
\newtheorem{lemma}[proposition]{Lemma}

\newtheorem{example}{Example}

\newtheorem{remark}{Remark}

\newcommand{\be}{\begin{equation}}
\newcommand{\ee}{\end{equation}}
\newcommand{\ba}{\begin{eqnarray}}
\newcommand{\ea}{\end{eqnarray}}
\newcommand{\bas}{\begin{eqnarray*}}
\newcommand{\eas}{\end{eqnarray*}}

\graphicspath{ {./Figures/} }

\def\T{{\mathcal T}}
\def\Re{{\mathbb R}}

\def\calL{{\mathcal L}}

\def\lmds{\texttt{lmds}}
\def\tlmds{\texttt{tlmds}}

\def\test{\texttt{test}}
\def\gap{\texttt{gap}}

\def\bfq{{\bf q}}

\def\bfs{{\bf s}}

\def\bfa{{\bf a}}
\def\bfb{{\bf b}}

\def\bfd{{\bf d}}

\def\bfx{{\bf x}}
\def\bfy{{\bf y}}
\def\bfz{{\bf z}}
\def\bfw{{\bf w}}
\def\bfu{{\bf u}}
\def\bfdelta{{\bf {\boldsymbol \delta}}}
\def\bfepsilon{{\bf {\boldsymbol \epsilon}}}
\def\bfzero{{\bf 0}}
\def\bfone{{\bf 1}}

\def\diag{\mbox{diag}}

\def\oa{{\overline{\bfa}}}

\begin{document}

\title{Landmark MDS Revisited for Sensor Network Localization with Anchors\thanks{This version: \today. This work was supported by the Hong Kong RGC General Research Fund (PolyU/15303124), and PolyU/P0044200, P0045347, and by
111 Project of China (B16002) and the National Science Foundation of China (12071022).}
  }

\author{Ting Ouyang,\thanks{T. Ouyang is with
		Department of Data Science and Artificial Intelligence,
		The Hong Kong Polytechnic University,
		 Hong Kong SAR, China. E-mail: tting.ouyang@connect.polyu.hk.}
\ \ Lingchen Kong\thanks{L. Kong is with School of Mathematics and Statistics, Beijing Jiaotong University, Beijing, 100044, China. E-mail: lchkong@bjtu.edu.cn.}
\ \ and \ \
Houduo Qi\thanks{H. Qi is with Department of Data Science and Artificial Intelligence, and Department of Applied Mathematics, The Hong Kong Polytechnic University, Hong Kong SAR, China. E-mail: {houduo.qi@polyu.edu.hk.}}
}



\maketitle

\begin{abstract}
The landmark multi-dimensional scaling (LMDS) is a leading method that embeds new points to an existing coordinate system based on observed distance information.
It has long been known as a variant of Nystr\"{o}m algorithm.
It was recently revealed that LMDS is Gower's method proposed in 1960s.
However, the relationship with other range-based embedding methods including
the least-squares (LS) has been unexplored,
proposing the question of which method to use in practice.
This paper provides a fresh look at those methods and explicitly
differentiates them in terms of the objectives they try to minimize.
For the first time for the case of single source localization,
we show that both LMDS and LS are generated from a same family of objectives, which balance
between length and angle
preservation among the embedding points.
Despite being nonconvex, the new objectives can be globally optimized
through a trust-region method.
An important result is that the LS solution can be thought as a regularized solution of LMDS.
Extension to the case of multiple source localization is also explored.
Comprehensive numerical results demonstrate the quality of the proposed objectives and the
efficiency of the trust-region method.
\end{abstract}


\section{Introduction}
\IEEEPARstart{M}{apping} new points to an existing coordinate system is a fundamental problem in
many applications such as data visualization \cite{buja2008data}, clustering \cite{little2023analysis}, 
localization \cite{biswas2006semidefinite}, and the internet of things \cite{ng2002predicting}, to just name a few. 
This problem is often cast as a generalized Procrustes problem \cite{gower1975generalized, gower2004procrustes} and there exist many algorithms depending on what information is available to use.
We are interested in those distance-based methods, which make use of
the measured (noisy) distances from the new points to anchors in the coordinate system.
A leading method is the Landmark Multi-Dimensional Scaling (LMDS) proposed in \cite{de2004sparse}, which has long been known as a variant of Nystr\"{o}m algorithm  \cite{platt2005fastmap}.
Until recently \cite{delicado2024multidimensional},
It was revealed to be Gower's method \cite{gower1968adding} proposed in 1960s.
Another popular method is the least-squares (LS) based \cite{beck2008exact}.
Despite their extensive use in various applications, 
the relationship between LMDS and LS remains unknown. 
The lack of understanding of their respective purposes and comparison 
causes the practical issue as to when to use them. 
The purpose of the paper is to address this issue
by showing that they can be cast as special instances of a unified framework,
leading to new models and efficient algorithms.
Significant implications to other methods are also investigated.

\subsection{Literature review and motivations}

\subsubsection{The problem of SSL}
To facilitate our discussion, let us use the framework of sensor network localization 
of one unknown source for demonstration, 
also known as the Single Source Localization (SSL) problem.
Suppose there are $m$ known sensors in $r$-dimensional space
(e.g., $r=2$ or $3$) and their
coordinates are known: $\bfx_i \in \Re^r$, $i \in [m]$, where $[m]$ denotes the set of
$\{1, \ldots, m\}$. 
Those sensors are often called anchors \cite{biswas2006semidefinite} or landmarks
\cite{de2004sparse}.
There is one unknown source $\bfx \in \Re^r$, whose distances to the anchors can
be observed or calculated through other means.
We denote the noisy distance from $\bfx$ to $\bfx_i$ by $\delta_i$. 
We assume that it is an approximation to the squared Euclidean distance: $\delta_i \approx \| \bfx - \bfx_i \|^2$,
where $\| \cdot\|$ is the Euclidean norm in $\Re^r$.
The task is to localize the true position of $\bfx$ based on the anchor positions
$\bfx_i$ and the observed distances $\delta_i$, $i\in [m]$.

\subsubsection{Least-squares model}
The least-squares (LS) criterion \cite{beck2008exact} is often used for estimating $\bfx$:
\be \label{LSQ-Original}
\min_{\bfx \in \Re^r} \ f_{ls}(\bfx) := \frac 1{2m} \sum_{i=1}^m \Big(  
\| \bfx - \bfx_i \|^2 - \delta_i
\Big)^2 .
\ee 
We note that the errors between the squared distances were used above. This
leads to a differentiable objective, though nonconvex.
The LS objective is also known as the squared stress in \cite{borg2005modern}
and it has been widely used in various localization problems \cite{beck2008exact}.
Gradient descent algorithms and more advanced semi-definite programming
relaxation methods can be developed
\cite{biswas2006semidefinite, so2007theory,shi2023facial}.

\subsubsection{MDS-based models}
MDS offers a completely different approach.
It starts with the Euclidean distance matrix $D$ consisting of the squared pairwise
distances between the anchors: $D_{ij} = \| \bfx_i - \bfx_j \|^2$, $i, j \in [m]$.
It then generates a set of new embedding points $\bfa_i \in \Re^r$.
Basic theory of MDS \cite{cox2000multidimensional, borg2005modern}
ensures that the distances are kept: $\| \bfa_i - \bfa_j \|^2 = D_{ij}$ for
all $i, j \in [m]$.
In other words, MDS creates a new coordinates system decided by those $\bfa_i$.
It then maps a new point to this system based on the observed distance measurements 
$\delta_i$, $i\in [m]$. 
The map can be nonlinear in $\delta_i$ (e.g., ISOMAP in \cite{tenenbaum2000global}).
The map by LMDS \cite{de2004sparse} is linear and hence can be fast calculated.
Many new points can be mapped to
the new system at a low cost.
We will explain how the LMDS mapping is obtained. In particular, 
we will explain what information is
kept and what is lost?
If one wants to get corresponding embedding to the original system defined by $\bfx_i$, 
a Procrustes mapping is needed \cite{gower2004procrustes, qi2013lagrangian, kong2019classical}
and it will be reviewed  in the next section.

As long been known, Gower has tackled this problem in 1960s \cite{gower1968adding}.
His basic tool is through distance triangulation.
Gower's starting point is
to assume that the observed distances $\delta_i$ are accurate and Euclidean.  
Then those distances can be reproduced by a new
point $[\bfx; x_0]$ in the $(r+1)$-dimensional space to the augmented landmark points
$\oa_i := [\bfa_i; 0] \in \Re^{r+1}$
(one more dimension than the MDS coordinate system $\bfa_i$). 
Then the unknown $\bfx$ must satisfy a system of linear equations
and this leads to the Gower mapping. 
It was recently proved that Gower's mapping and LMDS mapping are the same 
\cite{delicado2024multidimensional}.

Strictly speaking, neither Gower's method nor LMDS follows the spirit of MDS, which
is to use a low-rank Gram matrix to approximate MDS similarities among all points (known and unknown). 
As rightly pointed out in \cite{platt2005fastmap}, LMDS is actually a Nystr\"{o}m algorithm, which only approximates 
one part of MDS similarity matrix. 
In contrast, Trosset and Priebe \cite{trosset2008out} follows the spirit of MDS
in approximating the whole MDS similarity matrix
to derive a mapping that is a solution of a nonconvex optimization problem.
To differentiate this key difference, we refer to the method of Trosset and Priebe
as Total LMDS (T-LMDS).
Throwing away the nonconvex part \cite{trosset2008out}, it becomes 
a method by Anderson and Robinson 
\cite{anderson2003generalized}, referred to as the AR method in this paper. 

On the surface, LS model and the MDS-based models look quite distinct.
In fact, existing research on them went different ways.
We will show that they can be generated from a unified framework.
The numerical implication is significant as we summarize below.

\subsection{Main contributions}

We make three major contributions. 
The first one is about proposing a unified framework that cast all
the models reviewed above as special cases.
Our new perspective is to set every model in the MDS coordinate system
$\{\bfa_i\}$. 
The main findings are summarized in Table~\ref{Table-Objectives}.
In particular, we will show that the least-squares method \eqref{LSQ-Original} stated
in the original coordinate system determined by the anchors $\{ \bfx_i\}_{i \in [m]}$
can be equivalently stated in the new system. 
This new characterization, which is a special case of
the weighted T-LMDS with $w= 2/m$, allows us to extend 
existing numerical methods for the LS model to the weighted T-LMDS.
The weighted T-LMDS has a weight $w>0$ and it is a diagonally weighted version of
Trosset and Priebe \cite{trosset2008out} (T-LMDS uses $w=1$).
The weight reflects the importance of the measurements in $\delta_i$.
The quantities $b_0$ and $b_i$, $i\in [m]$ in Table-\ref{Table-Objectives} will be defined later.

\begin{table}[ht]
	\caption{Objectives of the Embedding Methods}	
	\centering
	\renewcommand{\arraystretch}{1.5}
	\begin{tabular}{l|l}
		\hline	
		Embedding Method & Objective to Minimize \\ [0.5ex]
		\hline
		LMDS/Gower/AR & $\sum_{i=1}^m \left(
		\langle \bfa_i,\; \bfx \rangle - b_i
		\right)^2 $ \\ \hline
		LS  & $\frac 12 (\| \bfx\|^2 - b_0)^2 + 
		\frac 2m \sum_{i=1}^m  (\langle \bfa_i, \; \bfx \rangle - b_i)^2$ \\ \hline
		T-LMDS  & $\frac 12 (\| \bfx\|^2 - b_0)^2 + 
		\sum_{i=1}^m  (\langle \bfa_i, \; \bfx \rangle - b_i)^2$ \\ \hline
		Weighted T-LMDS & $\frac 12 (\| \bfx\|^2 - b_0)^2 + 
		w \sum_{i=1}^m  (\langle \bfa_i, \; \bfx \rangle - b_i)^2$ \\ \hline
	\end{tabular}
	\label{Table-Objectives}
\end{table}

\noindent
It can be seen that each objective has two parts. One concerns the length of the new
embedding point $\| \bfx\|$.
The other part concerns the cosine similarity of the new point $\bfx$ with the existing
points $\bfa_i$. 
It can be easily seen that LMDS sets the weights $w_i$ to infinity so that the length
preservation part completely vanishes. 

Our second contribution is on the extension of the trust region method developed for 
the LS model \eqref{LSQ-Original} in \cite{beck2008exact} to the weighted T-LMDS model.
A major result in \cite{beck2008exact} is that the nonconvex LS model \eqref{LSQ-Original}
can be globally solved via a trust-region method.
The well-known ``hard case'' in the method was left out and was numerically verified 
that it rarely happened.
In this paper, we prove that the hard case can also be solved to its global solution, 
thanks to the nice properties of the MDS coordinate system $\{\bfa_i\}$. 
Since this is a major claim, we provide a water-tight mathematical proof. 

The third contribution is on extending T-LMDS to the case of
multiple source localization, which allows missing observations.
We show that the developed trust region method can be repeatedly applied to solve this 
case. Numerical experiments on both synthetic and real data demonstrate the 
quality of the proposed model and the efficiency of the algorithm.


\subsection{Organization} 

The paper is organized as follows. 
In Sect.~\ref{Section-Review}, we review the 
MDS-based embedding methods in detail and lay foundation for further development.
We propose the unified framework in Sect.~\ref{Section-T-LMDS}, where we prove that
the LS model \eqref{LSQ-Original} can be reformulated as a special case.
In Sect.~\ref{Section-Bisection}, we show the trust region method (TRM) can solve the
weighted T-LMDS model to its global optimality.
In particular, the hard case in the TRM can be efficiently solved.
In Sect.~\ref{Section-Multiple}, we extend the T-LMDS to 
the case of multiple source localization, which permits missing observations.
We further show that the TRM can be repeatedly applied to this case.
Finally, we report our comprehensive numerical experiments in
Sect.~\ref{Section-Numerical} in comparison with a few leading solvers.
We conclude the paper in Sect.~\ref{Section-Conclusion}.

\section{Review of MDS Embedding Methods} \label{Section-Review}

Let us recall the data we have.
There are $m$ landmarks $\bfx_i \in \Re^r$, $i\in [m]$ and an unknown source $\bfx$.
The distance measurement from $\bfx$ to the landmarks $\bfx_i$ is $\delta_i$:
$\delta_i \approx \| \bfx_i - \bfx\|^2$.
We let
\[
  \bfdelta^\top := [\delta_1, \ldots, \delta_m]
  \quad \mbox{and} \quad
  D = \Big(
    \| \bfx_i - \bfx_j \|^2
  \Big)_{i,j=1}^m,
\]
where $\bfdelta^\top$ is the transpose of $\bfdelta$ ($\bfdelta$ is a column vector).
This section reviews four MDS-based methods that locate $\bfx$.
They are the classical MDS \cite{borg2005modern}, LMDS \cite{de2004sparse},
Gower's method \cite{gower1968adding}, and T-LMDS \cite{trosset2008out}.
We also show that the method in \cite{anderson2003generalized} is actually
Gower's method.

\subsection{MDS representation of landmarks} \label{Subsection-MDS}

\subsubsection{Generating the new points $\{\bfa_i\}$}

The basic theory of MDS \cite{cox2000multidimensional, borg2005modern} says that
a new set of points $\bfa_i$ can be generated from the distance matrix $D$ to preserve
the known distances. It can be done as follows. 
Let
\[
  B := -\frac 12 J_m D J_m \quad \mbox{with}
  \ \ J_m := I_m - \frac 1m \bfone_m \bfone_m^\top ,
\]
where ''$:=$'' means ''define'', and $I_m$ is the $m \times m$ identity matrix and
$\bfone_m$ is the column vector of all ones of length $m$.
It is known that the matrix $B$ (also known as MDS similarity matrix) is positive semidefinite  and its rank is $r$.
The spectral decomposition of  $B$ is given by
\begin{align*}
  B &= \begin{bmatrix}
  	\bfu_1, \cdots, \bfu_r
  \end{bmatrix} \begin{bmatrix}
  \lambda_1 & & \\
   & \ddots & \\
   & & \lambda_r
  \end{bmatrix} \begin{bmatrix}
  \bfu_1^\top \\
  \vdots \\
  \bfu_r^\top
  \end{bmatrix}\\
  &= U_r \Lambda U_r^\top ,
\end{align*}
where $\lambda_1 \ge \ldots \lambda_r >0$ are the positive eigenvalues of $B$ and
their corresponding normalized eigenvectors are $\bfu_i$, $i=1, \ldots, r$.
For simplicity, we denote $\Lambda$ be the diagonal matrix consisting of the positive 
eigenvalues $\lambda_i$, $i=1, \ldots, r$ and $U_r$ consisting of corresponding eigenvectors $\bfu_i$.
The MDS embedding $\bfa_i \in \Re^r$ of the landmarks are given by
\be \label{MDS-Embedding}
 A := \begin{bmatrix}
 	\bfa_1, \cdots, \bfa_m
 \end{bmatrix} 
 = \Lambda^{1/2} U_r^\top .
\ee 
The most important properties of the new embedding are the distance-preserving and the centralization conditions:
\be \label{Centralization}
\begin{aligned}
    \| \bfa_i - \bfa_j \|^2 &= \| \bfx_i - \bfx_j \|^2 = D_{ij}, \ \
 i, j=1, \ldots, m, \\
  A \bfone_m &=  \bfa_1 + \cdots + \bfa_m = 0 .
\end{aligned}
\ee
This is because $\bfone_m$ is the eigenvector of $B$ corresponding to its zero eigenvalue. Consequently, $\bfu_i^\top \bfone_m = 0$ for all $i=1, \cdots, m$ and $A \bfone_m =0$.
The set of $\{ \bfa_i\}_{i=1}^m$ is often called the MDS coordinate system \cite{gower1966some}.

\subsubsection{The length of $\bfa_i$}

The length of each $\bfa_i$ can be directly calculated through the elements in $D$
without needing to actually calculate its coordinates.  Note that
\[
  \| \bfa_i - \bfa_j \|^2 = \| \bfa_i \|^2 + \| \bfa_j \|^2 - 2 \langle \bfa_i, \; \bfa_j \rangle .
\]
Summing over the index $j$ to get
\be \label{Eq-A}
\begin{aligned}
  \sum_{j=1}^m  \| \bfa_i - \bfa_j \|^2 &= m \| \bfa_i \|^2 + 
  \sum_{j=1}^m  \| \bfa_j \|^2 - 2 \langle \bfa_i, \; \sum_{j=1}^m \bfa_j \rangle\\
  &\stackrel{\eqref{Centralization}}{=}
  m \| \bfa_i \|^2 + 
  \sum_{j=1}^m  \| \bfa_j \|^2 .
  \end{aligned}
\ee 
Summing the above equation over the index $i$ to get
\[
 \sum_{i=1}^m \sum_{j=1}^m  \| \bfa_i - \bfa_j \|^2
 = 2m \sum_{j=1}^m  \| \bfa_j \|^2 .
\]
Consequently, we have
\[
  \sum_{j=1}^m  \| \bfa_j \|^2 = \frac 1{2m} \sum_{i=1}^m \sum_{j=1}^m  \| \bfa_i - \bfa_j \|^2
  \stackrel{\eqref{Centralization}}{=} \frac 1{2m} \bfone_m^\top D \bfone_m .
\]
It follows from \eqref{Eq-A} that
\begin{eqnarray} \label{Length-a}
 \| \bfa_i \|^2 &=& \frac 1m \sum_{j=1}^m  \| \bfa_i - \bfa_j \|^2 
 - \frac 1m \sum_{j=1}^m  \|  \bfa_j \|^2 \nonumber \\ 
 &=& \frac 1m \Big( D \bfone_m \Big)_i - 
 \frac 1{2m^2} (\bfone_m^\top D \bfone_m) .
\end{eqnarray}

\subsubsection{Mapping MDS anchors $\{\bfa_i \}$ to the original anchors $\{ \bfx_i\}$}

We note that the computation of the length of each $\bfa_i$ only used the relationship
\eqref{Centralization}. We now define
\[
 \bfx_0 := (\bfx_1 + \cdots + \bfx_m)/m 
 \quad \mbox{and} \quad 
 \widetilde{\bfx}_i := \bfx_i - \bfx_0, \ i \in [m].
\]
Then the set of points $\{ \widetilde{\bfx}_i\}$ also satisfy the
identities in \eqref{Centralization}. 
We must have $\| \bfa_i \| = \| \widetilde{\bfx}_i\|$ for $i \in [m]$. 
We define the following matrices and do a singular-value decomposition:
\[
  \widetilde{X} := [\widetilde{\bfx}_1, \ldots, \widetilde{\bfx}_m], \quad
  A \widetilde{X}^\top = U \Sigma V^\top , \quad
  P := V U^\top,
\]
where $U, V$ are $r \times r$ orthogonal matrices and $\Sigma$ is the $r \times r$ diagonal matrix containing the singular values of $A \widetilde{X}^\top$.
According to \cite[Prop.~4.1]{qi2013lagrangian},  we must have
\be \label{Mapping-x}
 \bfx_i = \T(\bfa_i) := P \bfa_i + \bfx_0, \quad i \in [m]
\ee 
where the linear mapping $\T(\cdot)$ will be used to map any new point in the MDS coordinate system $\{ \bfa_i\}_{i=1}^m$ to the anchor coordinate system 
$\{\bfx_i\}_{i=1}^m$.
\subsection{LMDS}

Having described MDS above, it is straightforward to describe LMDS.
Let
\[
 \bfdelta_0 := \frac 1m D \bfone_m
  \quad \mbox{and} \quad
 \calL^{\#}_r := \begin{bmatrix}
 	\frac 1{\sqrt{\lambda_1}} \bfu_1^\top \\
 	\vdots\\
 	\frac 1{\sqrt{\lambda_r}} \bfu_r^\top
 \end{bmatrix}.
\]
In fact, $\bfdelta_0$ is the column mean of the landmark distance matrix $D$
and $\calL^{\#}_r$ is the generalized inverse of the MDS embedding matrix
$A$ in \eqref{MDS-Embedding}. 
The embedding point in the MDS coordinate system 
for the unknown source by the landmark MDS \cite{de2004sparse}
is 
\[ 
  \hat{\bfx}_\lmds := \calL_\lmds (D, \bfdelta) 
  = -\frac 12 \calL^{\#}_r (\bfdelta - \bfdelta_0) .
\] 
Consequently, $\bfx_\lmds = \T(\hat{\bfx}_\lmds)$ is the embedding point
of the unknown source in the anchor system.
\subsection{The Gower embedding}

We omit the technical development presented by Gower \cite{gower1968adding}
and just list the equation that the new point $\bfx$ must satisfy :
\[ 
  A^\top \bfx = \frac 12 J_m ( \bfd - \bfdelta), \quad \mbox{with} \ \ 
  \bfd := \Big(\| \bfa_i\|^2 \Big)_{i=1}^m ,
\] 
where $\bfd$ is the column vector of squared length of each point $\bfa_i$.
Multiplying the above equation from the left by $A$ yields
\be \label{Gower-Normal-Eq}
   A A^\top \bfx = \frac 12 A J_m (\bfd - \bfdelta).
\ee 
Using the facts
\be \label{AJm}
   A J_m = A - \frac 1m A\bfone_m \bfone_m^\top \stackrel{\eqref{Centralization}}{=}  A
\ee
and
\be \label{Eq-AAT}
A A^\top \stackrel{\eqref{MDS-Embedding}}{=}
\Lambda^{1/2} U_r^\top U_r \Lambda^{1/2}
= \Lambda ,
\ee 
Equation \eqref{Gower-Normal-Eq} yields the solution of Gower embedding:
\[ 
 \hat{\bfx}_G = \calL_G(D, \bfdelta) := - \frac 12
\Lambda^{-1} A(\bfdelta - \bfd) .
\] 
With the notation $\bfdelta_0$,
the identity in \eqref{Length-a} implies
$$
\bfd = \bfdelta_0 - (1/2m^2) (\bfone_m^\top D \bfone_m) \bfone_m .
$$
Using $A\bfone_m = 0$, we get
\begin{small}
    \begin{align*}
	\hat{\bfx}_G &= \calL_G (D_0, \bfdelta)  \\
	&= -\frac 12 \Lambda^{-1} A \left(  \Big( \bfdelta - \bfdelta_0 + (1/2m^2) (\bfone_m^\top D \bfone_m)\bfone_m \right) \\
	&\stackrel{\eqref{Centralization}}{=} -\frac 12 \Lambda^{-1} A \Big(
	 \bfdelta - \bfdelta_0
	\Big) \\
	&= -\frac 12 \calL^{\#}_r (\bfdelta - \bfdelta_0) \qquad \mbox{(used $\Lambda^{-1} A = \calL^{\#}_r$ )}\\
	&= \calL_\lmds (D, \bfdelta).
\end{align*}
\end{small}
That is, the Gower method and LMDS produce the same embedding point.
This result was recently proved in \cite{delicado2024multidimensional}.

\subsection{Trosset-Priebe embedding}

Although both Gower mapping and LMDS are closely related to MDS embedding, they are
not strictly following the spirit of MDS.
Trosset and Priebe \cite{trosset2008out} do. We present their method below.
Let
   \[
\Delta := \begin{bmatrix}
	D & \bfdelta \\
	\bfdelta^\top & 0
\end{bmatrix}, \ 
\bfs := \frac 1m \begin{bmatrix}
	\bfone_m \\ 0
\end{bmatrix} 
\ 
\mbox{and} \  
J_\bfs := I_{m+1} - \bfs \bfone_{m+1}^\top .
\] 

If $\bfdelta$ contains the true squared Euclidean distances from $\bfx$ to each
landmarks, then the matrix $\overline{B} := -(1/2) J_\bfs^\top \Delta J_\bfs$ must
be positive semidefinite, see \cite[Thm.~2]{gower1982euclidean}.
It is easy to calculate
\begin{small}
\be
\label{eqbb0}
 \overline{B} = \begin{bmatrix}
 	-\frac 12 J_m D J_m &  \bfb \\ 
 	\bfb^\top  &  b_0
 \end{bmatrix} 
 = \begin{bmatrix}
 	A^\top A & \bfb \\
 	\bfb^\top & b_0
 \end{bmatrix},
\ee 
\end{small}
where 
\[
 b_0 := \frac 1m \bfone_m^\top \bfdelta - \frac 1{2m^2} \bfone_m^\top D \bfone_m,\
 \bfb := \frac 1{2}  J_m ( \bfdelta_0 -  \bfdelta),
\]
and we have used the MDS embedding for the matrix $B = -(1/2) J_m D J_m =A^\top A$
described in Subsection~\ref{Subsection-MDS}.
Therefore, when $\bfdelta$ contains noises it is natural to approximate
$\overline{B}$ by the Gram matrix
\[
  \begin{bmatrix}
  	A^\top  \\
  	\bfx^\top 
  \end{bmatrix}
  \begin{bmatrix}
  	A & \bfx
  \end{bmatrix}.
\]
This leads to the following least-squares problem:
\[
 \min_\bfx \ \frac 12 \left\|
  \overline{B} - \begin{bmatrix}
  	A^\top \\
  	\bfx^\top 
  \end{bmatrix}
  \begin{bmatrix}
  	A & \bfx
  \end{bmatrix}
 \right\|^2 , 
\]
which is equivalent to
\be \label{TP-LSQ}
\min_\bfx \ f_{TP}(\bfx) :=\frac 12 \Big( \| \bfx\|^2 - b_0 \Big)^2 +  \sum_{i=1}^m \Big(
\langle \bfa_i, \bfx \rangle - b_i
\Big)^2 .
\ee 
Trosset and Priebe \cite{trosset2008out} then left the problem to be solved 
by an optimization solver. 
We note that the optimization problem is nonconvex (i.e., $4$th-order polynomial 
optimization).
A global solution is often hard to obtain.

Dropping the first term in \eqref{TP-LSQ}, we get the model by Anderson and Robinson \cite{anderson2003generalized}:
\[  
\min_\bfx \ f_{AR}(\bfx) :=  \sum_{i=1}^m \Big(
\langle \bfa_i, \bfx \rangle - b_i
\Big)^2 .
\] 
Its optimality condition is
$
 AA^\top \bfx = A \bfb,
$
which leads to its embedding point:
\begin{align*}
	\hat{\bfx}_{AR} & = \calL_{AR}(\bfx) 
	= (A A^\top)^{-1} A\bfb
	  \stackrel{\eqref{Eq-AAT}}{=} \Lambda^{-1} A \bfb\\
	&\stackrel{\eqref{AJm}}{=} \Lambda^{-1} A \Big( \bfdelta_0 - \bfdelta  \Big) / 2 \\
	& = -\frac 12 \calL^{\#}_r (\bfdelta - \bfdelta_0)\\ 
	&= \calL_{\lmds} (D, \bfdelta) = \calL_G (D, \bfdelta) .
\end{align*}

We have shown that the three existing embedding methods (Gower embedding, Anderson-Robinson embedding and LMDS) are all
the same.
Since the full information of the similarity matrix $\overline{B}$ (including
the similarity matrix $B$ from the landmarks) was used, 
we refer to \eqref{TP-LSQ} as Total LMDS (T-LMDS).

\section{Weighted T-LMDS} \label{Section-T-LMDS}

In this part, we first propose a weighted version of Trosset-Priebe mode, denoted as
the weighted T-LMDS.
We then prove that the LS model \eqref{LSQ-Original} can be represented as an instance
of the weighted T-LMDS in the MDS coordinate system. 
We also make some comments on the choice of weights.

\subsection{Weighted T-LMDS}

It is natural to give a weight $w>0$ to the unknown source $\bfx$ to reflect its measurement quality.
Let
$\bfw^\top := [\bfone_m^\top, w]$. Consider the weighted problem of Trosset and Priebe
\cite{trosset2008out}:
\[
\min_\bfx \ \frac 12 \left\| \diag(\bfw) \left(  
\overline{B} - \begin{bmatrix}
	A^\top \\
	\bfx^\top 
\end{bmatrix}
\begin{bmatrix}
	A & \bfx
\end{bmatrix} \right) \diag(\bfw) 
\right\|^2 .
\]
Using the identities in \eqref{eqbb0}, this reduces to the problem:
\[ 
\min_\bfx \ \frac {w^2}2 \Big( \| \bfx\|^2 - b_0 \Big)^2 + w \sum_{i=1}^m \Big(
\langle \bfa_i, \bfx \rangle - b_i
\Big)^2 .
\]
Dividing $w^2$ on the objective and renaming $1/w$ by $w$, we arrive at the following weighted problem:
\be \label{Weighted-TP}
\min_\bfx \ \underbrace{\frac 12 \Big( \| \bfx\|^2 - b_0 \Big)^2}_{\mbox{length preservation}} +  \underbrace{w \sum_{i=1}^m   \Big(
\langle \bfa_i, \bfx \rangle - b_i
\Big)^2}_{\mbox{angle preservation}} .
\ee 
For easy reference, we call \eqref{Weighted-TP} the weighted T-LMDS.
We will see below that the least-squares problem \eqref{LSQ-Original}
is a special case.

\subsection{Least-squares as a weighted T-LMDS}

It is surprising to relate the least-squares problem 
\eqref{LSQ-Original} to the weighted T-LMDS problem \eqref{Weighted-TP}
as they are constructed from different criteria. 
We state this reformulation as follows.

\begin{theorem} \label{Prop-LSQ-Weighted-MDS}
The least-squares problem \eqref{LSQ-Original} is equivalent to the following problem
\be \label{f_lsq}
\min_{\widetilde{\bfx}} \ f_{ls} (\widetilde{\bfx}) = \frac 12 (\| \widetilde{\bfx}\|^2 - b_0)^2 + 
\frac 2m \sum_{i=1}^m  (\langle \bfa_i, \; \widetilde{\bfx} \rangle - b_i)^2 .
\ee 
Furthermore, suppose $\widetilde{\bfx}^*$ is an optimal solution of \eqref{f_lsq}, then
$
  \bfx^* = \T(\widetilde{\bfx}^*) = P \widetilde{\bfx}^* + \bfx_0
$
is an optimal solution of \eqref{LSQ-Original}, where the mapping $P$ and $\bfx_0$ are defined in
\eqref{Mapping-x}.
\end{theorem}

{\bf Proof.}
We first make a variable transformation. Let
$
 \widetilde{\bfx} := P^\top ( \bfx - \bfx_0).
$ 
Using the fact $PP^\top =I_r$, we obtain
\begin{align*}
  \| \bfx_i - \bfx\|^2 &= \| P \bfa_i + \bfx_0 - \bfx \|^2 
 = \| P ( \bfa_i - P^\top  ( \bfx - \bfx_0) ) \|^2 \\
 & = \|  \bfa_i - \widetilde{\bfx}_i \|^2.
\end{align*} 
Therefore, the least-squares problem \eqref{LSQ-Original} becomes
\be \label{LSQ-in-A}
\min_{\widetilde{\bfx} \in \Re^r} \ f_{ls}(\widetilde{\bfx}) := \frac 1{2m} \sum_{i=1}^m \Big(  
\| \widetilde{\bfx} - \bfa_i \|^2 - \delta_i 
\Big)^2 .
\ee 
Using the identity \eqref{Length-a}, we get
\begin{align*}
	&\| \bfa_i\|^2 - \delta_i 
	= \frac 1m \Big( D \bfone_m \Big)_i - 
	\frac 1{2m^2} (\bfone_m^\top D \bfone_m) - \delta_i  \nonumber \\
	&= \underbrace{\left(
		\frac 1m \Big( D \bfone_m \Big)_i - \delta_i + 
		\frac 1m \bfone_m^\top \bfdelta -
		\frac 1{m^2} (\bfone_m^\top D \bfone_m)
		\right)}_{= 2b_i}\\
	&\quad - \underbrace{ \left(
		\frac 1m \bfone_m^\top \bfdelta - \frac 1{2m^2} (\bfone_m^\top D \bfone_m)
		\right)}_{= b_0} \nonumber \\
	&= 2b_i - b_0, 
\end{align*}
and
\[
\sum_{i=1}^m ( \delta_i - \| \bfa_i\|^2)
= mb_0 - 2 \bfone_m^\top \bfb = m b_0,
\]
where the last equation used the property $\bfone_m^\top \bfb =0$ due to the definition 
of $\bfb$.
For simplicity, let $y = \| \widetilde{\bfx} \|^2$, we further have
\begin{align*}
	& \sum_{i=1}^m  \Big(  
	\| \widetilde{\bfx} - \bfa_i \|^2 - \delta_i 
	\Big)^2  \\
	& = \sum_{i=1}^m ( \| \widetilde{\bfx} \|^2 + \| \bfa_i\|^2 - \delta_i 
	- 2 \langle \bfa_i, \; \widetilde{\bfx} \rangle  )^2 \\
	& = m y^2 - 2y \sum_{i=1}^m ( \delta_i - \| \bfa_i\|^2)
	+ \sum_{i=1}^m ( \delta_i - \| \bfa_i\|^2)^2 
	 \\
	&- \underbrace{4y  \langle \sum_{i=1}^m \bfa_i, \; \widetilde{\bfx} \rangle}_{=0 \ \mbox{by} \ \eqref{Centralization}}
	+ \sum_{i=1}^m 4 ( \delta_i - \| \bfa_i\|^2) \langle \bfa_i, \; \widetilde{\bfx} \rangle + \sum_{i=1}^m 4 \langle \bfa_i, \; \widetilde{\bfx} \rangle^2 \\
	&= my^2 - 2my b_0 + 4 \sum_{i=1}^m \left(
	\langle \bfa_i, \; \widetilde{\bfx} \rangle - \frac 12 (\| \bfa_i\|^2 - \delta_i)
	\right)^2 \\
	&= m (y^2 - 2yb_0) + 4 \sum_{i=1}^m \left(
	\langle \bfa_i, \; \widetilde{\bfx} \rangle - (b_i - \frac 12 b_0)
	\right)^2 \\
	&= m (y - b_0)^2 - mb_0^2 \\
    & + 4 \sum_{i=1}^m \left[ (\langle \bfa_i, \; \widetilde{\bfx} \rangle - b_i)^2 
	- b_0 (\langle \bfa_i, \; \widetilde{\bfx} \rangle - b_i ) + \frac 14 b_0^2
	\right] 
\end{align*}
\begin{align*}
	&= m (y - b_0)^2 - mb_0^2 + 4 \sum_{i=1}^m  (\langle \bfa_i, \; \widetilde{\bfx} \rangle - b_i)^2 + mb_0^2  \\
    &  - 4b_0 
	\underbrace{\left(
		\langle \sum_{i=1}^m \bfa_i, \; \widetilde{\bfx} \rangle - \bfone_m^\top \bfb
		\right)}_{=0 \ \mbox{by} \ \eqref{Centralization}\ \mbox{and} \ \bfone_m^\top \bfb=0}\\
	&= m (y - b_0)^2  + 4 \sum_{i=1}^m  (\langle \bfa_i, \; \widetilde{\bfx} \rangle - b_i)^2 .
\end{align*}
Therefore, the least-squares embedding problem \eqref{LSQ-in-A} can be
reformulated as Problem \eqref{f_lsq} after substituting $y = \| \widetilde{\bfx}\|^2$ back.
Let $\widetilde{\bfx}^*$ be an optimal solution of \eqref{f_lsq}. 
We map back to get $\bfx^* = P \widetilde{\bfx}^* + \bfx_0$ as an
optimal solution of the original problem \eqref{LSQ-Original}.
\hfill $\Box$ \\

What we have proved above is that the LS model \eqref{LSQ-Original}, when put in the MDS
coordinate system $\{\bfa_i\}_{i=1}^m$, can be represented as 
a weighted T-LMDS problem \eqref{Weighted-TP} with $w= 2/m$.
We will see that in this new system, the LS problem \eqref{f_lsq} is easier to
solve than its original version \eqref{LSQ-Original}.
We summarize the four models we considered so far in Table~\ref{Table-Objectives} in Introduction.
They are all represented in the MDS coordinate system.

\subsection{Choice of weights} \label{Section-Weight-choice}

The weighting parameter $w$ in \eqref{Weighted-TP} serves to balance two complementary geometric constraints: length preservation and angle preservation. 
By adjusting this weight, we can adaptively control the algorithm's sensitivity to noise in each constraint type, thereby improving accuracy.

We like to show that the noises in $\bfdelta$ affect the length and the angle terms differently. 
Let $\bfd^*$ be the column vector representing the true squared distances between the unknown point and the landmarks:
$d^*_i = \| \bfx - \bfx_i\|^2$, $i \in [m]$.
We assume the noises are in the additive model:
\[
  \bfdelta = \bfd^* + \bfepsilon,
\]
where $\bfepsilon \in \Re^m$ is the noise vector.
Let $\bfb^*$ and $b_0^*$ be obtained by replacing $\bfdelta$ in $\bfb$ and $b_0$ by $\bfd^*$.
That is, $\bfb$ and $b_0$ are the respective noisy observation of $\bfb^*$ and $b_0^*$. 
Then we have
\[ 
|b_0-b_0^*| = \frac 1 m |\bfone_m^{\top}(\bfdelta-\bfd^*)|.
\] 
Noise impacts $b 
_0$ primarily through the sum $\bfone_m^{\top}\bfepsilon$, which is the sum of noises.
The angle preservation term can be reformulated as follows:
	\[
	\begin{aligned}
		\sum\limits_{i=1}^m\big(\langle \bfa_i,\bfx\rangle-b_i\big)^2 & =  \|A^{\top}\bfx-\bfb\|^2\\
		& = \| A^{\top}\bfx\|^2-2\langle A\bfb,\bfx\rangle + \|\bfb\|^2.
	\end{aligned}
	\] 
Hence, noise primarily affects the angle preservation term through the product $A\bfb$ (while $\|\bfb\|^2$ contains noise, it remains a constant in optimization). 
Therefore, we can measure the noise level:
\[ 
2\|A\bfb-A\bfb^*\| = \|AJ_m(\bfd^*-\bfdelta)\| \stackrel{\eqref{Centralization}}{=}
\|A \bfepsilon\|.
\] 
The following two observations are obvious:
\begin{enumerate}[label=C.\arabic*]
	\item \label{condition1} $\bfone_m^{\top}\bfepsilon = 0$ has no impact on the length preservation term. 
	\item \label{condition2}$A\bfepsilon = 0$ has no impact on the angle preservation term.
\end{enumerate}
Of course, in practice we do not know whether the noise vector $\bfepsilon$ satisfies
the two conditions above.
But they suggest which term we should place more weight on if we know something about
the noise level.
In particular, we have the following comments on the four models.

\begin{remark} \label{Remark-on-noise}
(i) Since LMDS only focuses on the angle preservation part, its performance can be 
satisfactory when the noise in the angle part is small.
(ii) The least-squares model places less emphasis on the angle component compared 
with the T-LMDS model as the LS weight $2/m \le 1$ and $1$ is the weight used in T-LMDS.
(iii) The weighted T-LMDS has an adjustable weight on the angle-preservation part.
It will be more flexible under various kinds of circumstances.
\end{remark}

\section{Reformulation as Trust-Region Problem} \label{Section-Bisection}

While all models reviewed so far can be united in the form of the weighted problem
\eqref{Weighted-TP}, LMDS (i.e., Gower's method) has a closed-form formula.
But all others must rely on an iterative algorithm.
Beck et.al., \cite{beck2008exact} proposed a trust region method for the problem
\eqref{LSQ-Original}. 
Under certain regularity conditions, the standard trust region algorithm can be applied
with its ``hard case'' removed due to the regularity condition enforced.
In this part, we will show that the same trust region method can be applied to the weighted problem \eqref{Weighted-TP} and we do not need any regularity condition to remove the hard case. 
Since this is a strong claim, 
we give a detailed proof of how this can be done. 
We first recall two major results in the trust region method from \cite{more1993generalizations}.

\subsection{Two results on generalized trust-region problem}

Consider the generalized trust region problem:
\be \label{GTR}
 \min_{\bfz \in \Re^n} \ q(\bfz) \quad \mbox{subject to} \quad c(\bfz) = 0,
\ee 
where both $q(\cdot)$ and $c(\cdot)$ are quadratic functions:
\[
  q(\bfz) := \frac 12 \bfz^\top Q \bfz - \bfq^\top \bfz
  \quad \mbox{and} \quad
  c(\bfz) := \frac 12 \bfz^\top Q_0 \bfz - \bfq_0^\top \bfz + c_0,
\] 
where $Q_0, Q$ are $n \times n$ symmetric matrices, $\bfq_0, \bfq \in \Re^n$, and $c_0 \in \Re$.
The global solution of \eqref{GTR}  is characterized as follows.

\begin{lemma}\label{Lemma-More}
%
%
\cite[Thm.~3.2]{more1993generalizations} 
	(Characterization of global minimizer)
	Suppose $\nabla^2 c \not=0$ and
\[
  \inf\{ c(\bfz) \; | \ \bfz \in \Re^n \}
  < 0 < \sup\{ c(\bfz) \; | \ \bfz \in \Re^n \}.
\]
A vector $\bfz^*$ is a global minimizer of \eqref{GTR} if and only if 
there is a multiplier $\lambda^* \in \Re$ such that the Karush-Kuhn-Tucker
(KKT) 
conditions 
\be \label{KKT}
  \left\{ \begin{array}{rr}
  	\nabla q(\bfz^*) + \lambda^* \nabla c(\bfz^*) &= 0\\ [1ex]
  	\nabla^2 q(\bfz^*) + \lambda^* \nabla^2 c(\bfz^*)& \succeq 0\\ [1ex]
  	   c(\bfz^*) &=0
  \end{array}
  \right .
\ee
are satisfied, where a matrix $M \succeq 0$ means it is positive semidefinite.

\end{lemma}

In order to find $\bfz^*$, we let $\bfz(\lambda)$ be the solution of the
equation for a given $\lambda \in \Re$:
\be \label{KKT-1}
 \nabla q(\bfz(\lambda) ) + \lambda \nabla c(\bfz (\lambda)) = 0
\ee
and define
\[
 I_{PD} := \left\{
  \lambda \in \Re \ | \ Q + \lambda Q_0\  \mbox{is positive definite}
 \right\}.
\]
It is known that $I_{PD}$ is an open interval.
The following result is the basis of the efficient bisection method for finding the optimal $\lambda^*$ and $\bfz^*$.

\begin{lemma} \label{Lemma-Decreasing}
	\cite[Thm.~5.2]{more1993generalizations}
Let $\lambda \in I_{PD}$. Then $\bfz(\lambda)$ is well defined.
We further let
\[
  \phi(\lambda) := c (\bfz(\lambda)).
\]
Then $\phi(\cdot)$ is strictly decreasing on $I_{PD}$ unless $\bfz(\cdot)$ is constant on $I_{PD}$ with
\[
  \nabla q(\bfz(\lambda)) = 0, \quad \nabla c(\bfz(\lambda)) = 0
  \quad \mbox{for all} \ \lambda \in I_{PD}.
\]
	
\end{lemma}

\subsection{Application to the weighted problem} 

We now apply the results above to the weighted problem \eqref{Weighted-TP}.
Let
\[
 \overline{A} := \begin{bmatrix}
 	\sqrt{2w} A & \bfzero_r \\
 	\bfzero_m^\top & 1
 \end{bmatrix}, \  \overline{\bfb} := \begin{bmatrix}
 \sqrt{2w} \bfb \\
  b_0
 \end{bmatrix},
\]
and 
\[
Q_0 := \begin{bmatrix}
 	I_r & \\
 	 & 0
 \end{bmatrix}, \  \bfq_0 := \begin{bmatrix}
 \bfzero_m \\ 1/2
 \end{bmatrix}, \  \bfz := \begin{bmatrix}
 \bfx \\ y
 \end{bmatrix}.
\]
Then the weighted problem \eqref{Weighted-TP}
can be further cast as the generalized trust-region problem:
\be \label{TP-Trust}
\begin{array}{ll}
 \min_{\bfz}  &  \ q(\bfz) = \frac 12 \| \overline{A}^\top \bfz - \overline{\bfb} \|^2 \\ [0.6ex]
 \mbox{s.t.} & \ c(\bfz) = \frac 12 \bfz^\top Q_0 \bfz - \bfq_0^\top \bfz = 0.
\end{array} 
\ee 
Note that $c(\bfz) =0$ means $y = \| \bfx\|^2$.
In this setting, the derivatives of both functions take simple forms:
\begin{align*}
 Q:=  &	\nabla^2 q(\bfz) = \overline{A} \overline{A}^\top =
	\begin{bmatrix}
		(2w) A A^\top & \\
		  & 1 
	\end{bmatrix} = \begin{bmatrix}
	(2w) \Lambda & \\
	 & 1
	\end{bmatrix} \\
   &	\nabla q (\bfz) = \overline{A} \overline{A}^\top \bfz -
	\overline{A} \overline{\bfb}
	= \begin{bmatrix}
		(2w) \Lambda \bfx \\
		y 
	\end{bmatrix} -
	\begin{bmatrix}
		2w A \bfb \\
	    b_0
	\end{bmatrix} \\
	& \nabla^2 c(\bfz) = Q_0, \quad \nabla c(\bfz) = Q_0 \bfz - \bfq_0 =
	\begin{bmatrix}
		\bfx \\
		-\frac 12
	\end{bmatrix} .
\end{align*}
We have the following result.

\begin{theorem} \label{Prop-Weighted-Problem}
Problem \eqref{TP-Trust} has an optimal solution.
Furthermore, $\bfz^* = (\bfx^*, y^*)$ is its global solution if and only if
there exists a Lagrange multiplier $\lambda^* \in \Re$ such that
the following conditions hold:
\be \label{KKT-Weighted}
 \left\{
 \begin{array}{l}
  \Big(  
   \Lambda + \frac 1{2w} \lambda^* I_r 
  \Big) \bfx^* = A \bfb , \\ [0.8ex]
  y^* = b_0 + \frac 12 \lambda^*, \ \ \lambda^* \ge -2w \lambda_r, \\ [0.8ex]
  y^* = \| \bfx^*\|^2 .
 \end{array} 
 \right .
\ee 
	
\end{theorem}

{\bf Proof.}
The result is a restatement of Lemma~\ref{Lemma-More}.
Since the Hessian matrix $\nabla^2 q(\bfz)$ is positive definite,
the function $q(\bfz)$ is coercive (i.e., $q(\bfz) \rightarrow +\infty$ when
$\| \bfz\|\rightarrow +\infty$), and hence Problem \eqref{TP-Trust} must have
a global minimizer.
When $\bfx=0$ and $y>0$, we have $c(\bfz) = -y <0$ and when $\bfx=0$, $y <0$,
we have $c(\bfz) >0$. 
Moreover, $Q_0 = \nabla^2 c(\bfz) \not=0$.
Therefore, the conditions in Lemma~\ref{Lemma-More} are
satisfied. It is safe to use the KKT conditions \eqref{KKT} for \eqref{TP-Trust}.
They are stated as follows.
\[
 \left\{ \begin{array}{l}
 	 ( \Lambda + \frac 1{2w} \lambda^* I_r ) \bfx^* = A\bfb, \ 
 	 \bfy^* = b_0 + \frac 12 \lambda^* \\ [1ex]
 	 \begin{bmatrix}
 	 	(2w) \Lambda & \\
 	 	  & 1
 	 \end{bmatrix} + \lambda^* \begin{bmatrix}
 	 I_r & \\
 	 & 0 
 	 \end{bmatrix}  \succeq 0\\ [1ex]
 	 y^* = \| \bfx^*\|.
 \end{array} 
 \right .
\]
Simplifying the above system leads to 
\eqref{KKT-Weighted}. \hfill $\Box$\\

The task now is to find $(\lambda^*, \bfx^*, y^*)$ that satisfies \eqref{KKT-Weighted}.
We note that both $\bfx^*$ and $y^*$ depend on $\lambda^*$ and 
\eqref{KKT-Weighted} is nonlinear. Hence, it is usually hard to solve.
We also note that $I_{PD}$ takes the following form:
\[\begin{aligned}
	I_{PD} &= \left\{
	\lambda \in \Re \ | \ Q + \lambda Q_0 \succ 0
	\right\}\\
	&= \left\{
	\lambda \in \Re \ | \ (2w) \Lambda + \lambda I_r \succ 0
	\right\}\\ 
	&= (-(2w)\lambda_r, \infty),
\end{aligned}
\]
where $\lambda_r>0$ is the smallest eigenvalue in $\Lambda$.
For the choice of $\lambda^*$,
there are only two cases: $\lambda^* = -2w\lambda_r$ or $\lambda^* \in I_{PD}$.
The first case is known as the ``hard case'' in the trust region method, while the second
case is often referred to as ``easy case''. 
Fortunately, for the hard case, we can give it a complete characterization 
as we explain below.


\subsubsection{Hard case}
Suppose the optimal $\lambda^* = -2w \lambda_r$. 
In this case, we can find a necessary and sufficient condition that any solution $\bfx^*$ will have to satisfy.
Let $\lambda_k$ be the second smallest eigenvalue in $\lambda_1 \ge \cdots \ge \lambda_r$, which comes from the decomposition \eqref{MDS-Embedding}.
That is $\lambda_k$ is the smallest eigenvalue such that $\lambda_k > \lambda_r$
(the situation that all positive eigenvalues are equal cannot happen in the hard case unless
$\bfx^*=0$ is the optimal solution). 
Let $\Lambda_k$ be the diagonal matrix of the first $k$ positive eigenvalues $\lambda_1 \ge \ldots \ge \lambda_k$.
Let $\widehat{\bfb} := A\bfb$ and $\widehat{b}_{[1:k]}$ be the subvector of
$\bfb$ consisting of the first $k$ components of $\bfb$. 
Then the first equation of \eqref{KKT-Weighted} becomes
\[
  \Big( \Lambda_k - \lambda_r I_k \Big) \bfx^*_{[1:k]} = \widehat{\bfb}_{[1:k]}.
\]
We have $\bfx^*_{[1:k]} = ( \Lambda_k - \lambda_r I_k )^{-1} \widehat{\bfb}_{[1:k]}$.
The necessary condition that $\bfx^*$ must satisfy is as follows:
\be \label{HardCaseCondition}
 \left\{
 \begin{array}{l}
 	\bfx^*_{[1:k]} = ( \Lambda_k - \lambda_r I_k )^{-1} \widehat{\bfb}_{[1:k]} \\
 	\widehat{\bfb}_{[(k+1):r]} = 0 \\
 	y^* = b_0 - w \lambda_r > 0 \\
 	y^* \ge \| \bfx^*_{[1:k]} \|^2 .
 \end{array} 
 \right .
\ee  
This condition is also sufficient because we simply set 
 $x^*_{k+1} = \sqrt{ y^* - \| \bfx^*_{[1:k]}\|^2}$ and $x^*_i =0$ for 
$i > k+1$. The resulting $\bfx^*$ will satisfy the KKT condition \eqref{KKT-Weighted}.

The condition in \eqref{HardCaseCondition} makes the hard case unlikely to happen. 
Anyway,
we can quickly check whether it happens by computing a feasible solution $\bfx^*$ above.
We note that this hard case was omitted in \cite{beck2008exact}. 
Instead, they assume some regularity condition to ensure only the easy case happens.

\subsubsection{Easy case}
For this case, we have $\lambda^* > -(2w) \lambda_r$.
That is, the optimal $\lambda^*$ falls in the open interval $I_{PD}=(-2w \lambda_r, \infty)$.
We will search this $\lambda^*$.
For $\lambda \in I_{PD}$, we let $\bfz(\lambda) = (\bfx(\lambda), y(\lambda))$ be 
the solution of the following equations:
\[
( \Lambda + \frac 1{2w} \lambda I_r ) \bfx (\lambda) = A\bfb, \quad 
\bfy(\lambda) = b_0 + \frac 12 \lambda .
\]
The solution is given by
\[
 \bfx(\lambda) = \Big( \Lambda + \frac 1{2w} \lambda I_r
 \Big)^{-1}  A\bfb, \quad y(\lambda) = b_0 + \frac 12 \lambda.
\]
It is easy to see that $\nabla c (\bfz(\lambda)) \not=0$ for any $\lambda \in I_{PD}$.
Lemma~\ref{Lemma-Decreasing} implies that 
\be \label{phi}
 \psi(\lambda) := \phi(c(\bfz(\lambda)))
= \| \bfx(\lambda)\|^2 
  - \frac 12 \lambda - b_0 
\ee
is strictly decreasing in $\lambda$ over the interval $I_{PD}$. 
The optimal $\lambda^*$ satisfies
$\psi(\lambda^*) =0$. 
The strict monotonicity of $\psi(\lambda)$ suggests that 
the bisection method is efficient for searching for
this $\lambda^*$. For its use, we derive an upper bound for
$\lambda^*$.
Let
\[
  \lambda_{\max} := 2\max\left\{
    \left(  \frac 1{ \lambda_r + 1} \right)^2 \| A \bfb \|^2 - b_0,\; w
  \right\}.
\]
Then we obtain
\begin{align*}
   \psi(\lambda_{\max})	&= c(\bfz(\lambda_{\max}))\\
	&= \| \bfx( \lambda_{\max}) \|^2 - \frac 12 \lambda_{\max} - b_0 \\
	&= \sum_{i=1}^r \left( \frac{1}{\lambda_i + (1/(2w))\lambda_{\max}}  \right)^2 
	   (A \bfb)_i^2 - \lambda_{\max} - b_0 \\
	&\le  \left( \frac{1}{\lambda_r + 1} \right)^2 \sum_{i=1}^k (A \bfb)_i^2 - \frac 12 \lambda_{\max} - b_0 \\ 
	&=  \left( \frac{1}{\lambda_r + 1} \right)^2 \|A \bfb \|^2 - b_0 - \frac 12\lambda_{\max}   \le 0 . 
\end{align*}
Therefore, $\lambda^* \in I_* := (-(2w)\lambda_r, \lambda_{\max})$.
We may apply the bisection method over the interval $I_*$ to find the
optimal $\lambda^*$ so that $c(\bfx(\lambda^*))=0$. Once $\lambda^*$ is obtained, we get the optimal 
embedding point
\be \label{TLMDS-x}
  \hat{\bfx}_\tlmds = \calL_\tlmds(D, \bfdelta) := \Big( \Lambda + \frac 1{2w}\lambda^* I_r \Big)^{-1} A \bfb .
\ee
\begin{remark} \label{Remark-TLMDS}
(i) The computation of the matrix inverse in $\bfx(\lambda)$ is easy as the matrix is 
diagonal. Therefore, the bisection search is very fast.
(ii) It is interesting to compare the solution formula \eqref{TLMDS-x} with that of LMDS:
\[
 \calL_\lmds(D, \bfdelta) = \Lambda^{-1} A \bfb.
\]
$\calL_\tlmds(D, \bfdelta)$ can be treated as a regularized version of $\calL_\lmds$.
We note that $\lambda^*$ may be negative.
\end{remark}

Alg.~\ref{B-SSL} summarizes the Bisection method for the reformulated problem \eqref{TP-Trust}
of SSL, denoted as B-SSL($\overline{A}, \overline{\bfb}, \epsilon$),
where $\epsilon>0$ is the given tolerance for terminating the bisection search
for the output $\hat{\bfx}_\tlmds$.
Finally, we may map $\hat{\bfx}_\tlmds$ by \eqref{Mapping-x} to get the final embedding $\bfx_\tlmds$ in the anchor system.

\begin{algorithm}[ht]
	\caption{(B-SSL($\overline{A}, \overline{\bfb}, \epsilon$) for Problem \eqref{TP-Trust})}  \label{B-SSL}
	\begin{algorithmic}[1]
		\IF{(hard case): Conditions \eqref{HardCaseCondition} are satisfied} 
		\STATE Set $\bfx^*_{k+1}=\sqrt{y^*-\|\bfx^*_{[1:k]}\|^2}$ to get $\bfx^*$.
		\ELSE[(easy case):] 
		\STATE Set $\lambda_{\min} := -2w \lambda_r$.
		\STATE Compute
			$$
			\lambda_\test :=\frac 12 (\lambda_{\min} + \lambda_{\max}),
			\lambda_\gap := \frac 12 (\lambda_{\max} - \lambda_{\min}).
			$$	\label{step-bisec}
		\STATE If $\lambda_\gap > \epsilon$, update
			\[
			\left\{
			\begin{array}{ll}
				 \lambda_{\min} = \lambda_\test & \mbox{if} \ \psi(\lambda_\test) \ge 0 \\
				 \lambda_{\max} = \lambda_\test & \mbox{if} \ \psi(\lambda_\test) \le 0.
			\end{array} 
			\right.
			\]
		\STATE Go to step \ref{step-bisec}.
		\ENDIF	 
		
		\STATE Output: $\lambda^* = \lambda_\test$ and compute
		$\hat{\bfx}_\tlmds$ by \eqref{TLMDS-x}.
	\end{algorithmic}
\end{algorithm}


\section{Extension to Multiple Source Localization} \label{Section-Multiple}

This part deals with the case of multiple source localization (MSL)
 with range observations, which
may have missing values. We first set up the problem. We then describe how we solve it 
utilizing the method for the SSL.

\subsection{MSL: Multiple source localization} 
We use slightly different notation. As before, there are $m$ anchors, whose known 
positions are denoted as $\bfx_i \in \Re^r$.
We have $n >1$ unknown sources, whose positions are denoted as
$\bfy_i\in \Re^r$, $i \in [n]$.  Hence, there are a total of $N := m+n$ points in
the system. The matrix $D$ still denotes the Euclidean distance matrix from the
$m$ anchors. We have two more matrices: $m \times n$ matrix $E$ and $n \times n$ matrix 
$F$, which contain the observed values of the true distances from the anchors to the unknown sources, and among the unknown sources $Y:= [\bfy_1, \ldots, \bfy_n]$:
\[
E \approx \Big( \| \bfx_i - \bfy_j \|^2 \Big)_{i \in [m]}^{j \in [n]} \
\mbox{and} \ 
F \approx \Big( \| \bfy_i - \bfy_j \|^2 \Big)_{i \in [n]}^{j \in [n]}.
\]
Our basic assumption is that $E$ is fully observed (no missing observations from an
anchor to an unknown source).
But we allow $F$ to contain missing values. In particular, we define for each unknown resource $i \in [n]$
\[
 \Omega_i := \left\{
  j \in [n] \setminus  \{i\} \ | \ F_{ij} \ \mbox{is observed}
 \right\}.
\]
In other words, when $\Omega_i = \emptyset$, there were no communication between 
source $\bfy_i$ to other unknown sources. 
The problem is to find the locations of those unknown sources.

\subsection{Generalized T-LMDS}
We follow the spirit of T-LMDS for SSL to build it for MSL \cite{trosset2008out}.
Let
\[
  \Delta := \begin{bmatrix}
  	D & E \\
  	E^\top & F
  \end{bmatrix}, \quad \bfs := \begin{bmatrix}
   \frac 1m \bfone_m \\
   \bfzero_n
  \end{bmatrix} \quad \mbox{and} \quad
  J_\bfs := I_N - \bfs \bfone_{N}^\top,
\]
We obtain the approximation:
\be  \label{Block-Approximation}
	 \begin{bmatrix}
		A^\top \\
		Y^\top 
	\end{bmatrix} \begin{bmatrix}
		A, Y
	\end{bmatrix} 
	\approx -\frac 12 J_\bfs^\top \Delta J_\bfs 
	 = -\frac 12 \begin{bmatrix}
		J_m D J_m, & \overline{E} 	\\
		\overline{E}^\top,
		&  \overline{F} 
	\end{bmatrix},
\ee
where
\begin{align*}
	\overline{E} &:= J_m \Big( E - \frac 1m D \bfone_m \bfone_n^\top  \Big)\\
	\overline{F} &:= F - \frac 1m \Big(   
	E^\top \bfone_m \bfone_n^\top + \bfone_n \bfone_m^\top E
	\Big) + \frac{ \bfone_m^\top D \bfone_m }{m^2} 
	\bfone_n \bfone_n^\top .
\end{align*}
It is important to note that the matrix $\overline{E}$ is well defined because it
only depends on $E$ and $D$ (both are available).
Furthermore, the missing values say $F_{ij}$ in $F$ do not translate to other elements
of $\overline{F}$. That is, $\overline{F}_{ij}$ is available if and only if
$F_{ij}$ is available. 
In particular, $\overline{F}_{ii}$ is well defined because $F_{ii}=0$.
This key observation makes it possible that we only approximate the available data
in \eqref{Block-Approximation}.
The optimization problem from \eqref{Block-Approximation} takes the following form:
\be \label{MDS-Approximation}
  \min_Y \ \frac 12 \left\|
  -\frac 12 J_\bfs^\top \Delta J_\bfs - \begin{bmatrix}
  	A^\top \\
  	Y^\top 
  \end{bmatrix}
  \begin{bmatrix}
  	A & Y
  \end{bmatrix}
  \right\|^2 , 
\ee
If we only focus on the available data and use the identity 
$A^\top A = -J_mDJ_m/2$, the above problem reduces to
\begin{align} \label{Optimization-MSL}
& \min_Y \  \| A^\top Y - \overline{E} \|^2 \nonumber   \\
&  +   \frac 12 \sum_{i=1}^n 
          \left[   ( \|\bfy_i \|^2 - \overline{F}_{ii}) ^2 +
           \sum_{j \in \Omega_i}
 ( \langle \bfy_i, \bfy_j \rangle - \overline{F}_{ij} )^2 \right] .
\end{align}
We further let 
\[
E := \begin{bmatrix}
	\bfdelta_1, \cdots, \bfdelta_n
  \end{bmatrix},  \quad \bfb_i := \frac 12 J_m ( \bfdelta_0 - \bfdelta_i), \ i=1, \ldots, n.
\]
Then problem \eqref{Optimization-MSL} can be rewritten as
\be \label{MSL-Problem}
 \min_{Y} \; \sum_{i=1}^n \left[
   f_i(\bfy_i) + F_i(\bfy_i, Y_{(-i)})
 \right],
\ee 
where
\begin{align*}
	f_i(\bfy_i) &= \frac 12 ( \| \bfy_i \|^2 - \overline{F}_{ii})^2 + \| A^\top \bfy_i - \bfb_i \|^2 \\
	F_i(\bfy_i, Y_{(-i)}) &= \sum_{j \in \Omega_i} 
	( \langle \bfy_i, \bfy_j \rangle - \overline{F}_{ij} )^2 ,
\end{align*}
and $Y_{(-i)}$ is the matrix $Y$ with the $i$th column $\bfy_i$ removed.
This notation will be convenient for describing our alternating minimization algorithm 
below.

\subsection{B-MSL: Alternating bisection method}

Unlike the case of SSL, Problem \eqref{MSL-Problem} is more challenging and
is hard to compute its global solution. Fortunately, it has a good structure and
an alternating minimization algorithm can be applied as follows: Given the current
iterate $\widehat{Y}^k$, compute 
\begin{align*}
  \hat\bfy^{(k+1)}_i &= \arg\min_{\bfy} f_i^{(k)}(\bfy) := f_i(\bfy) + F_i(\bfy, Y^k_{(-i)}) .
\end{align*} 
This subproblem is a generalized trust region problem and can be efficiently solved
by the bisection method developed with suitable modifications. 
To see this, let
$$
\overline{A}_{i}^{k} := \begin{bmatrix}
    \sqrt{2} A & \sqrt{2}(\bfy_j^k)_{j \in \Omega_i} &\bfzero_r\\
    0 &  0 &1
\end{bmatrix},
\  \overline{\bfb}_{i} = \begin{bmatrix}
    \sqrt{2} \bfb^{(i)} \\
    \sqrt{2} ( \overline{F}_{ij} )_{j \in \Omega_i} \\
    \overline{F}_{ii}
\end{bmatrix},\ 
  $$
and
$$
Q_0 := \begin{bmatrix}
 	I_r & \\
 	 & 0
 \end{bmatrix}, \  \bfq_0 := \begin{bmatrix}
 \bfzero_m \\ 1/2
 \end{bmatrix}, \  \bfz_i := \begin{bmatrix}
 \bfy_i \\ y
 \end{bmatrix}.
$$
The subproblem can be written as 
\be \label{MSL-yi}
\begin{array}{ll}
\min_{\bfz_i} & \frac 1 2\left\|\overline{A}_{i}^{k}\bfz_i-\overline{\bfb}_{i}\right\|^2 \\
\mbox{s.t.} & c(\bfz_i) = \frac 12 \bfz_i^\top Q_0 \bfz_i - \bfq_0^\top \bfz_i = 0.
\end{array} 
\ee 
The resulting method B-MSL (Bisection method for MSL) is summarized as follows.
\begin{algorithm}[ht]
  \caption{B-MSL Algorithm}\label{B-MSL}
  \begin{algorithmic}[1]
    \REQUIRE Matrix $\Delta$, anchors $\bfx_i,i\in[m]$, tolerance $\epsilon>0$.
    \STATE \textbf{Initialization: }set $k=0$, $\widehat{Y}^{0} = \mathcal{L}_{\mathrm{lmds}}(D,E)$.
    \REPEAT
      \FOR{$i=1$ to $n$}
        \STATE use B-SSL($\overline{A}^k_i, \overline{\bfb}_i, \epsilon$)
        for \eqref{MSL-yi}
         to obtain $\hat\bfy_i^{\,k+1}$.
      \ENDFOR
      \STATE set $\widehat{Y}^{k+1} = [\hat\bfy_1^{\,k+1},\dots,\hat\bfy_n^{\,k+1}]$.
      \STATE $k \leftarrow k+1$.
    \UNTIL{$\|\widehat{Y}^{k}-\widehat{Y}^{k-1}\| \le \epsilon$}.
    \STATE\textbf{Mapping: } Map each point in $\widehat{Y}^{k}$ by \eqref{Mapping-x} to get $Y_{\texttt{tlmds}}$.
  \end{algorithmic}
\end{algorithm}

\section{Numerical Experiments} \label{Section-Numerical} 

This section presents extensive numerical experiments that evaluate the proposed method B-SSL for single source localization, and B-MSL for multiple source localization. All experiments were performed using MATLAB (R2024b) on a desktop computer equipped with an Apple M3 CPU and 16 GB of RAM.

\subsection{Comparison on single source localization}
In this section, we focus on testing the B-SSL algorithm for the single source localization problem. We compare it against the previously mentioned methods (LMDS/Gower/AR, T-LMDS, and LS (Least-squares) and also include the representative state-of-the-art Trilateration method  (Tri. for short) (Algorithm 1) \cite{larsson2025single}. 

\begin{example} \label{ex1}
	This example is from \cite{beck2008exact}. 
	There are five anchors on a plane ($r = 2$) and their positions are given by
	\[
	[\bfx_1,\bfx_2,\bfx_3,\bfx_4,\bfx_5] = \begin{bmatrix}
		-5&-12&-1&-9&-3\\
		-13&1&-5&-12&-12
	\end{bmatrix}.
	\]
	The position of the sensor is $\bfx = (-5,11)^{\top}$. We use various noise type to demonstrate the effect of weight in B-SSL. The noisy range measurement is given by $\bfdelta = \bfd^*+\bfepsilon$, where $\bfd^*$ represents the squared true sensor-anchor distances and $\bfepsilon$ is the noise term satisfying specific conditions.
	Especially, we use three different types of noise in this example:
	\begin{enumerate}
		\item[(i)] $\bfepsilon^{\top}\bfone_5 = 0$, $\bfepsilon$ satisfies condition \ref{condition1}.
		\item[(ii)] $A\bfepsilon = 0$, $\bfepsilon$ satisfies condition \ref{condition2}.
		\item[(iii)] $\epsilon_i\in N(0,\sigma^2)$, $\epsilon_i$ is randomly generated from some Gaussian distribution with zero mean and $\sigma$ standard deviation.
	\end{enumerate}
	For types (i) and (ii), the noise vector $\epsilon$ is generated by orthogonally projecting a random Gaussian vector onto the subspace satisfying the required condition. All results are averaged over 500 Monte Carlo trials. Table \ref{tab:ssl} reports the average localization errors and the average weight $\overline{w}$ for B-SSL. We omit CPU times from the table, as the differences among the five algorithms are negligible. For cases (i) and (iii), the weight $w$ of B-SSL is given by a grid search over the interval $[0.01,\ 2]$ with step size 0.1. For case (ii), we simply set the weight to a sufficiently large value (e.g., $w =  10^4$).
	
\end{example}
\begin{table}[htbp]
	\centering
	\caption{Averaged position error of the five algorithms on Example \ref{ex1} with different noise type and level and average weight $\overline{w}$ for B-SSL}
	\label{tab:ssl}
	\setlength{\tabcolsep}{3.5pt} 
	\small 
	\begin{tabular}{cllllll}
		\hline
		$\sigma$ & B-SSL  & T-LMDS & LS & LMDS & Tri & $\overline{w}$ \\ \hline
		& \multicolumn{5}{c}{Noise Type (i)} & \\
		0.01 & \textbf{4.22E-04} & 4.43E-04 & 4.38E-04 & 6.35E-04 & 4.61E-04 & 0.87 \\
		0.1  & \textbf{4.29E-03} & 4.53E-03 & 4.45E-03 & 6.76E-03 & 4.58E-03 & 0.81 \\
		1    & \textbf{4.47E-02} & 4.71E-02 & 4.63E-02 & 6.90E-02 & 4.80E-02 & 0.78 \\ \hline
		& \multicolumn{5}{c}{Noise Type (ii)} & \\
		0.01 & \textbf{4.81E-05} & 1.05E-04 & 1.13E-04 & \textbf{4.81E-05} & 1.64E-04 & 10000 \\
		0.1  & \textbf{5.11E-04} & 1.05E-03 & 1.12E-03 & \textbf{5.11E-04} & 1.67E-03 & 10000 \\
		1    & \textbf{5.27E-03} & 1.02E-02 & 1.09E-02 & \textbf{5.27E-03} & 1.75E-02 & 10000 \\ \hline
		& \multicolumn{5}{c}{Noise Type (iii)} & \\
		0.01 & \textbf{4.45E-04} & 4.66E-04 & 4.65E-04 & 6.35E-04 & 4.96E-04 & 0.92 \\
		0.1  & \textbf{4.45E-03} & 4.68E-03 & 4.65E-03 & 6.76E-03 & 4.89E-03 & 0.87 \\
		1    & \textbf{4.64E-02} & 4.88E-02 & 4.85E-02 & 6.90E-02 & 5.09E-02 & 0.85 \\ \hline
	\end{tabular}
\end{table}

Table \ref{tab:ssl} shows that—when its weight is chosen appropriately—B-SSL consistently delivers the most accurate results among the five algorithms. By contrast, T-LMDS, LS, LMDS and Tri. each exhibit specific weaknesses when handling particular noise types, which reflects the comments in Remark \ref{Remark-on-noise}.
We make some further observations below.

{\bf On effects of the noise type}.
Under noise type (i), where length information is exact, B-SSL attains the highest accuracy. LS performs slightly better than T-LMDS. This is because LS assigns a weight of $w=2/m=0.4$ to the (noisy) angle-preservation term, whereas T-LMDS places a weight of $w=1$ on the angle term. Tri. is less accurate than the first three methods, and LMDS performs the poorest because it omits the length-preservation component and therefore cannot exploit the exact length information.

Under noise type (ii), where angle information is exact, the ordering changes: LMDS consistently achieves the best accuracy, and B-SSL can reach comparable performance by adjusting its weight. The remaining methods (T-LMDS, LS, and Tri.) degrade markedly.
Among them, T-LMDS fares best because it relies more heavily on angle information than LS, and Tri. ranks last.

Under noise type (iii), with purely random noise, B-SSL again attains higher accuracy than the alternatives when its weight is tuned appropriately. T-LMDS and LS outperform Tri., while LMDS yields the poorest results.



\subsection{Comparison on multiple source localization}
\begin{figure}[htbp]
	\centering
	\includegraphics[width=1\linewidth]{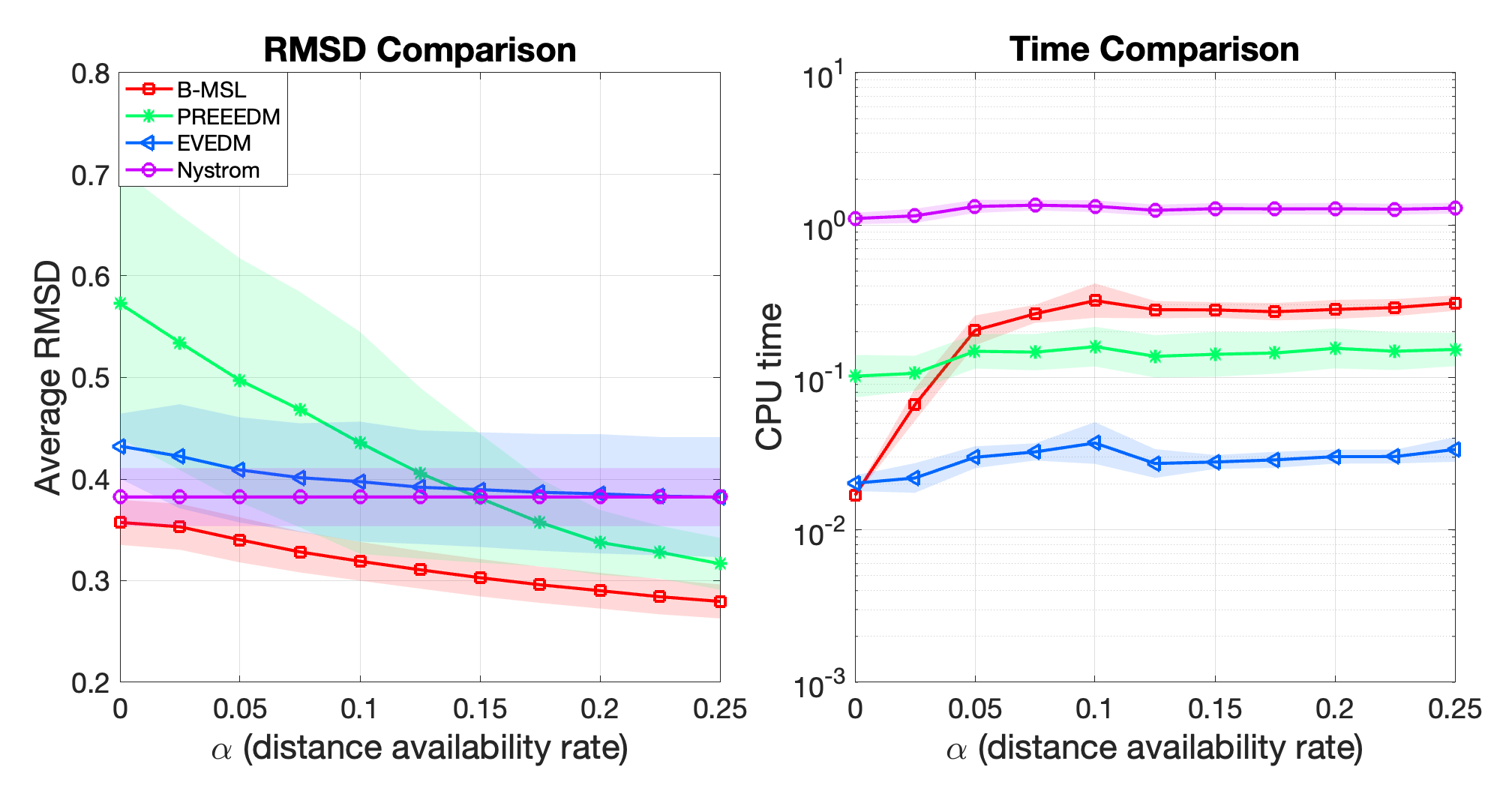}
	\caption{Average results for Example \ref{ex:synthetic} with $m=20$, $n=200$, $\sigma = 0.5$}
	\label{fig:ratedemo}
\end{figure}

\begin{figure}[htbp]
	\centering
	\includegraphics[width=1\linewidth]{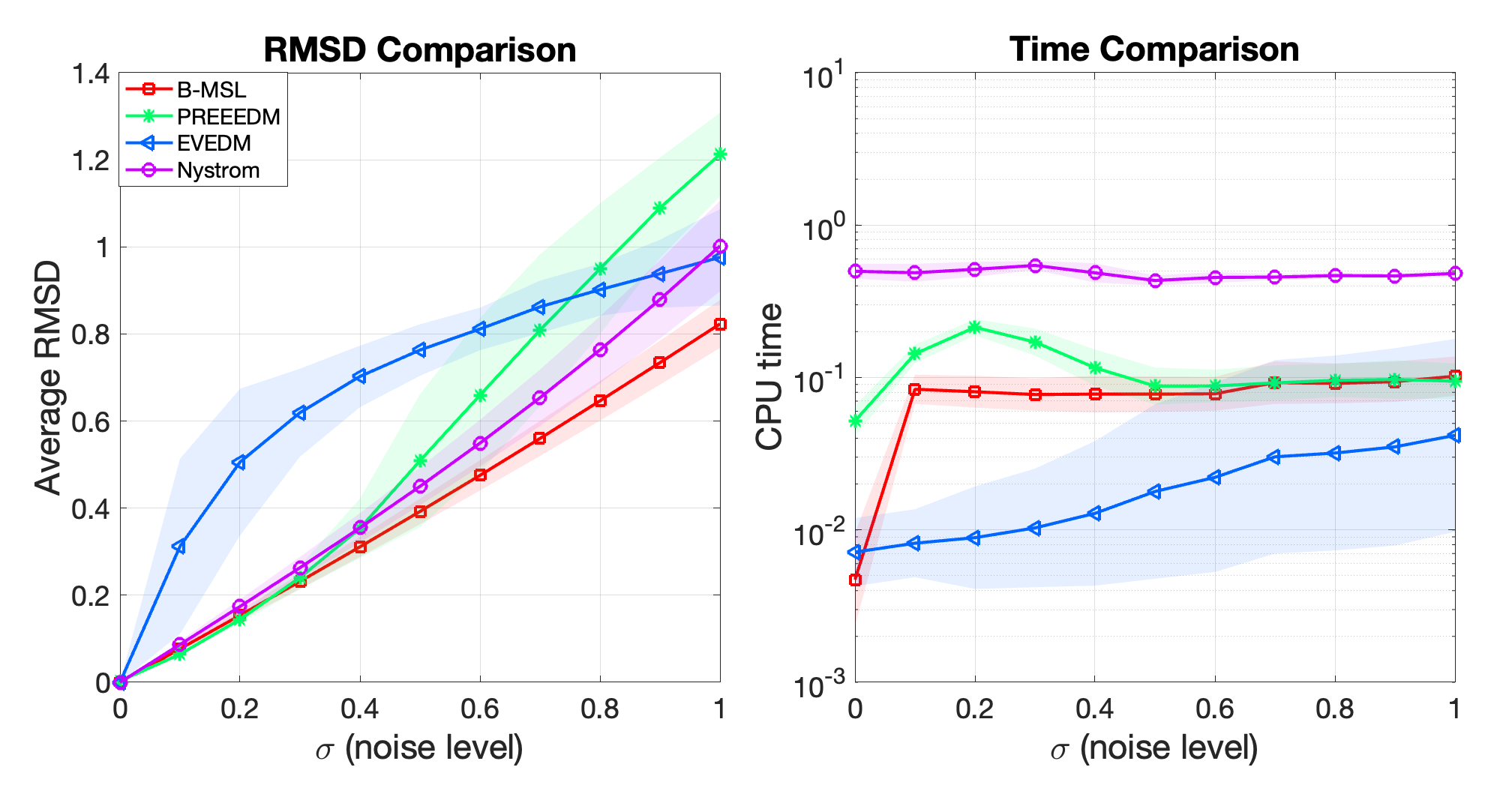}
	\caption{Average results for Example \ref{ex:synthetic} with $m=15$, $n=150$, $\alpha = 0.05$}
	\label{fig:noisedemo}
\end{figure}

\begin{figure*}[htbp]
	\centering
	\includegraphics[width=0.9\linewidth]
	{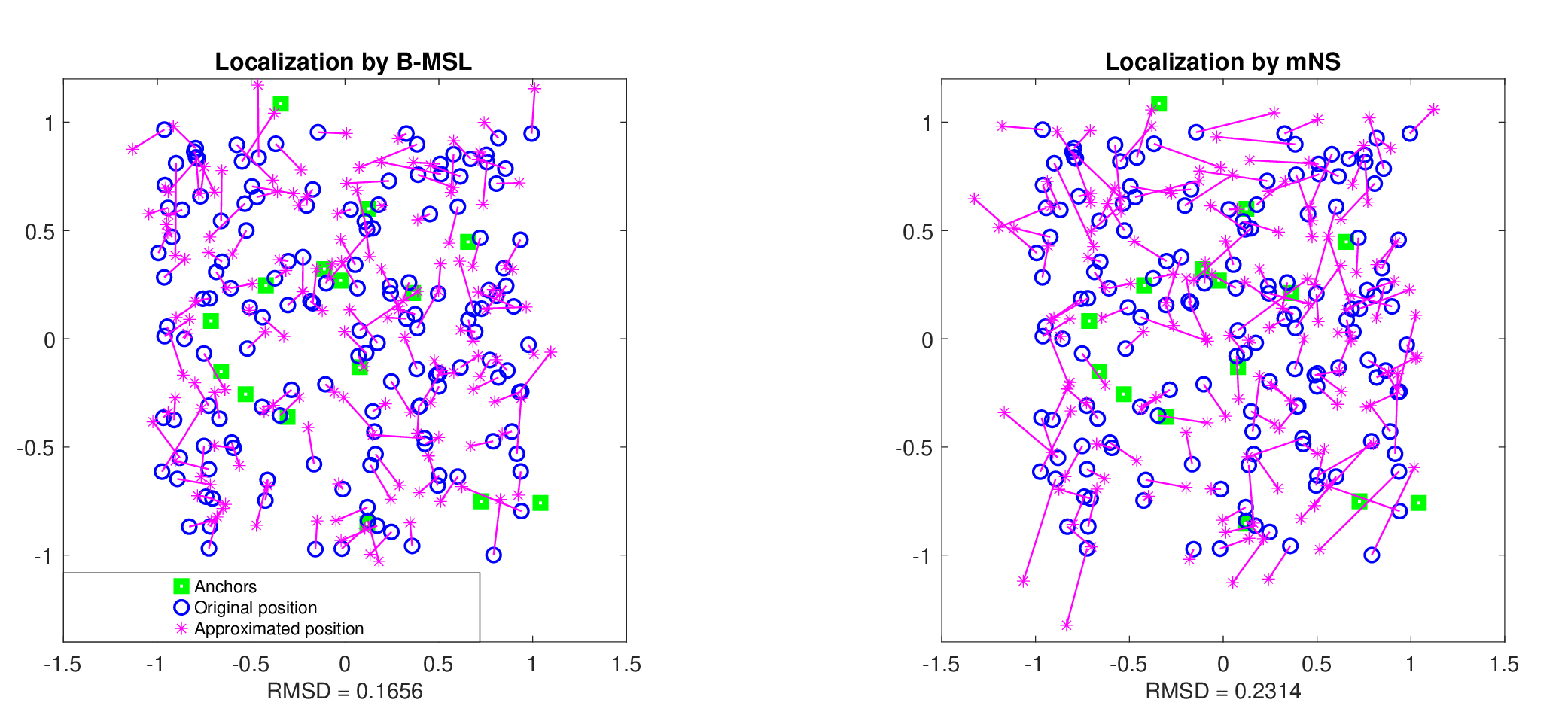}
	\caption{Localization results for B-MSL and mNS on Example \ref{ex:synthetic} with $m=15$, $n=150$, $\sigma = 0.25$, $\alpha = 0.15$}
	\label{fig:position}
\end{figure*}

In this section, we test and compare the B-MSL algorithm against other 3 state-of-the-art algorithms on both synthetic and real-world examples. The chosen algorithms are PREEEDM from \cite{zhou2020robust}, EVEDM (short for EepVecEDM)\cite{drusvyatskiy2017noisy}, and the modified Nystr\"om sampling method (mNS) from \cite{lichtenberg2024localization}.
Especially, we use Algorithm 1 in \cite{lichtenberg2024localization} since the anchor-anchor distances are all given. We report the computational time and the RMSD (Root of the Mean Squared Deviation), which is defined as $$\mbox{RMSD} = \sqrt{\frac{\|Y-Y_M\|^2}{n}},$$ where $Y$ is the true position matrix of the sensors and $Y_M$ is the approximation by different method.
\begin{example}\label{ex:synthetic}(Synthetic data) 
	In this example, we generate random instances based on an additive noise model: 
	\[
	\begin{cases}
		E_{ij} = (\|\bfx_i-\bfy_j\|+\sigma\epsilon_{ij})^2,&i\in[m],j\in[n],\\
		F_{ij}=(\|\bfy_i-\bfy_j\|+\sigma\epsilon_{ij})^2,&(i,j)\in\Omega.
	\end{cases}
	\]
	The sensors and the anchors are randomly distributed in the area $[-1,1]^2$. We evaluate the performance of all algorithms under three scenarios, with each result averaged over 100 independent runs:
	\begin{enumerate}
		\item Increasing the distance availability rate $\alpha$ from 0 to 0.25 (where $\alpha = |\Omega|/n^2$ represents the percentage of given sensor-to-sensor distances). Performance across all algorithms as $\alpha$ increases is shown in Fig. \ref{fig:ratedemo}.
		\item Increasing the noise level $\sigma$ from $0$ to $1$ in increments of $0.1$. Results are presented in Fig. \ref{fig:noisedemo}.
	\end{enumerate}
	The solid lines represent the mean RMSD (or CPU time), while the shaded regions indicate the standard deviation from the mean. As a specific example, Fig.  \ref{fig:position} plots the localizations achieved by B-MSL and mNS in a representative case.
\end{example}

\paragraph{Effects of distance availability rate $\alpha$} Fig. \ref{fig:ratedemo} shows that B-MSL, PREEEDM and EVEDM effectively leverage the available sensor–sensor distances. The RMSD for these methods decreases as the rate $\alpha$ increases, indicating improved accuracy.  In contrast, the RMSD of mNS changes little, as it does not incorporate the distance information between sources.   In terms of accuracy, B-MSL consistently achieves the highest accuracy across all rate. EVEDM and PREEEDM are the fastest methods, though they are less accurate than B-MSL. The Nystrom method is outperformed in terms of accuracy also with a slowest speed.

\paragraph{Effects of noise level $\sigma$} From Fig. \ref{fig:noisedemo} we can see that the accuracy of all tested algorithms dropped as the noise level $\sigma$ rise. However, B-MSL constantly demonstrates the highest accuracy most of the time (especially when the noise level is high), while other algorithms are less accurate under the same noise level. Also, B-MSL has a competitive speed compared with other algorithms. EVEDM is the fastest algorithm in terms of CPU time, followed by B-MSL and PREEEEDM, while mNS is the slowest.

From Fig. \ref{fig:position} we can see that B-MSL provides a more accurate reconstruction of the sensor positions compared to mNS, as indicated by the lower RMSD value. The lines in the figure, which represent the error for each sensor, are also visibly shorter for B-MSL, indicating its estimates are closer to the true positions. In contrast, mNS presents larger deviations for most sensors.

\begin{figure}[htbp]
	\centering
	\subfloat[B-MSL: RMSD = 0.1284 Å, Time = 1.0408 s]{%
		\includegraphics[width=0.22\textwidth]{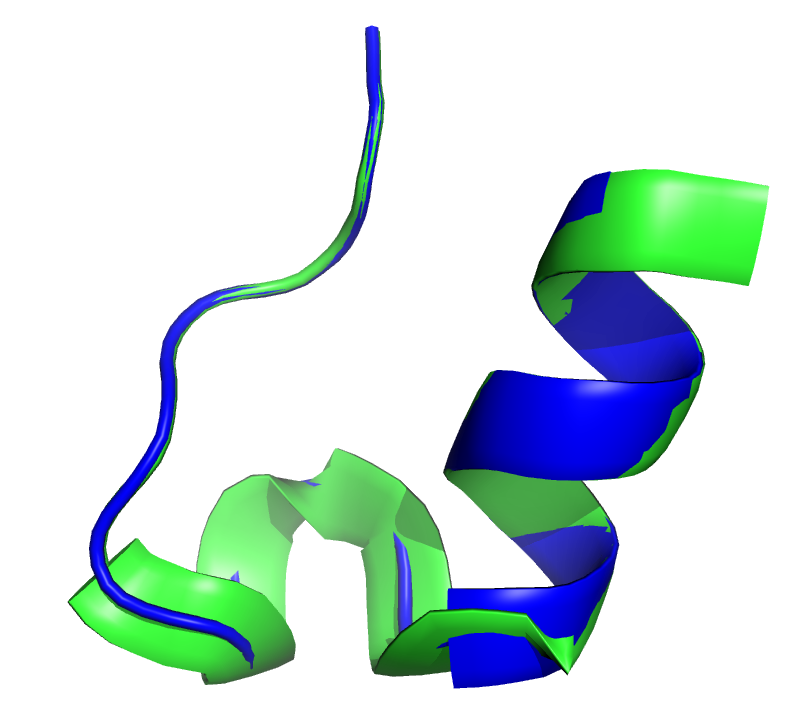}\label{fig:subfig7}}
	\hfill
	\subfloat[PREEEDM: RMSD = 0.9311 Å, Time = 0.0785 s]{%
		\includegraphics[width=0.22\textwidth]{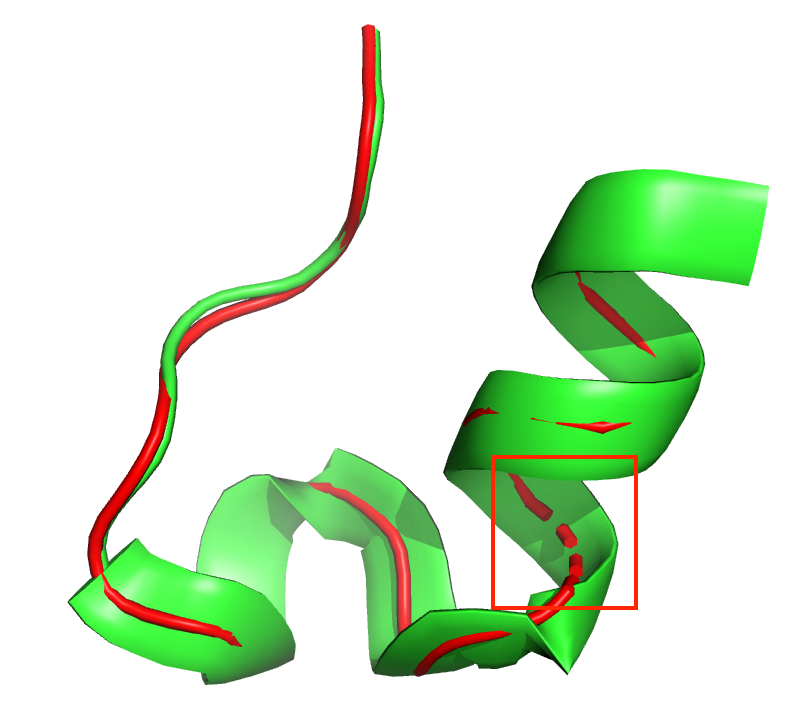}\label{fig:subfig8}}
	\vspace{0.5em}
	\subfloat[EVEDM: RMSD = 0.2516 Å, Time = 0.0070 s]{%
		\includegraphics[width=0.22\textwidth]{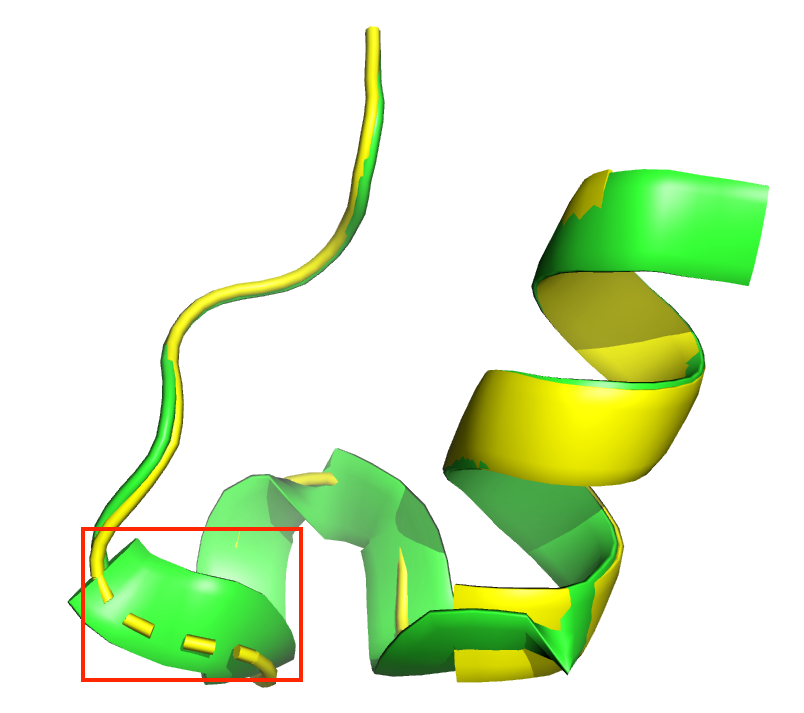}\label{fig:subfig9}}
	\hfill
	\subfloat[mNS: RMSD = 0.2197 Å, Time = 0.2017 s]{%
		\includegraphics[width=0.22\textwidth]{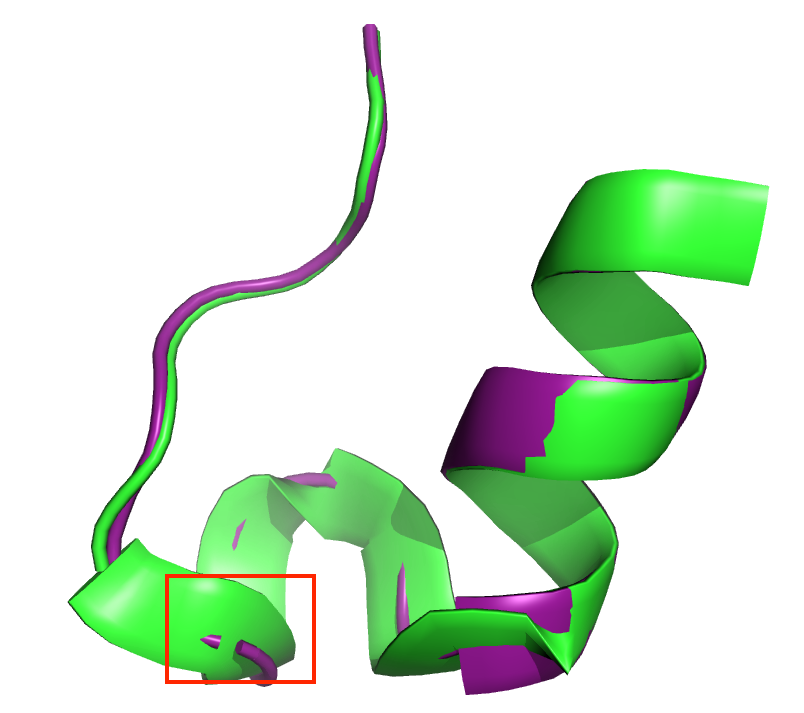}\label{fig:subfig10}}
	
	\caption{The structure of 1L2Y (in green) and the approximated structure by other algorithms (in blue, red, yellow and purple) with $\alpha=0.15$, $\sigma = 0.01$}
	\label{fig:protein}
\end{figure}

\begin{example}\label{ex:real}(Real data) In this example, we evaluate the algorithms on the representative protein data identified as 1L2Y \cite{Neidigh2002TrpCage, PDB1L2Y}. This protein contains 154 atoms (excluding most hydrogens). We select the first 5 atoms to be the anchors, and the others are the sources. The noises are added based on the same model in Example \ref{ex:synthetic}. The results are averaged over 100 runs and are recorded in Table \ref{tab:protein}. Especially, we visualize the example of $\alpha = 0.15$ and $\sigma = 0.01$ by PyMOL Molecular Graphics System \cite{pymol} in Fig. \ref{fig:protein}.
	\begin{table}[htbp]
		\caption{Results on Example \ref{ex:real} with $\alpha=0.1$\label{tab:protein}}
		\begin{tabular}{clllll}
			\hline
			$\sigma$    &      & B-MSL             & PREEEDM           & EVEDM    & mNS      \\ \hline
			0.01 & RMSD & \textbf{1.92E-01} & 3.34E+00          & 2.63E-01 & 2.16E-01 \\
			& Time & 1.16E+00          & 1.79E-01          & 6.67E-03 & 2.01E-01 \\
			0.1  & RMSD & \textbf{1.72E+00} & 3.65E+00          & 4.54E+00 & 2.16E+00 \\
			& Time & 9.67E-01          & 1.73E-01          & 1.13E-01 & 1.87E-01 \\
			1    & RMSD & 7.08E+00          & \textbf{3.56E+00} & 6.99E+00 & 2.16E+01 \\
			& Time & 9.33E-01          & 4.24E-01          & 8.08E-02 & 2.06E-01 \\ \hline
		\end{tabular}
	\end{table}
\end{example}

From Table \ref{tab:protein} we can see B-MSL significantly outperforms the others in accuracy under low to moderate noise levels ($\sigma=0.01,0.1$), though its computational time is generally not the fastest. When the noise level becomes high ($\sigma=1$), PREEEDM achieves the best accuracy. The mNS algorithm serves as a consistently fast and reliable option for low-to-moderate noise, where it secures the second-best accuracy, but its performance collapses dramatically in the high-noise scenario. EVEDM is the fastest algorithm under all cases, but provides relatively inaccurate results.

The results in Fig. \ref{fig:protein} show that B-MSL achieves the highest accuracy (lowest RMSD) but requires the longest computation time. Its reconstructed configurations closely match the true structure, with no dashed bonds. The other methods are significantly faster but exhibit lower accuracy and introduce structural discontinuities, as highlighted in the red boxes, which also means they are less accurate. PREEEDM is faster but less accurate, exhibiting noticeable deviations. EVEDM offers a reasonable trade-off between speed and accuracy, though it remains less precise than B-MSL. The mNS method attains moderate accuracy and is slower than both EVEDM and PREEEDM.

\section{Conclusion} \label{Section-Conclusion}

This paper reports a surprising connection between the landmark MDS-based methods and
the least-squares method for sensor network localization with anchors.
Casting in the classical MDS generated coordinates system for the anchors, all of those methods can be
generated from a class of weighted functions, which balance 
between the length preservation and the angle preservation from
the observed data.
This also provides a new perspective in understanding how noise was handled in
those methods. 
Some pay a higher weight to the length preservation, while the others place
a higher weight to the angle preservation.
The proposed weighted T-LMDS provides a flexible and simple way to adjust the weight.
 
Furthermore, the use of the MDS coordinate system also leads to a fast algorithm
that is based on bisections search. 
We proved that the bisection search is capable of finding the global solution for the
case of single unknown source despite it being a nonconvex optimization problem.
Extension to multiple source localization was also investigated through alternating
minimization. 
Numerical experiments
demonstrate that the new algorithms are effective for both cases.
A challenging theoretical question is whether the multiple source location problem
can be solved to its global optimality, even for the case $n=2$ (two unknowns).
We leave this for our future research.
 
\bibliographystyle{ieeetr}
\bibliography{mdsRefs}

\begin{thebibliography}{10}

\bibitem{buja2008data}
A.~Buja, D.~F. Swayne, M.~L. Littman, N.~Dean, H.~Hofmann, and L.~Chen, ``Data visualization with multidimensional scaling,'' {\em Journal of computational and graphical statistics}, vol.~17, no.~2, pp.~444--472, 2008.

\bibitem{little2023analysis}
A.~Little, Y.~Xie, and Q.~Sun, ``An analysis of classical multidimensional scaling with applications to clustering,'' {\em Information and Inference: A Journal of the IMA}, vol.~12, no.~1, pp.~72--112, 2023.

\bibitem{biswas2006semidefinite}
P.~Biswas, T.-C. Liang, K.-C. Toh, Y.~Ye, and T.-C. Wang, ``Semidefinite programming approaches for sensor network localization with noisy distance measurements,'' {\em IEEE transactions on automation science and engineering}, vol.~3, no.~4, pp.~360--371, 2006.

\bibitem{ng2002predicting}
T.~E. Ng and H.~Zhang, ``Predicting internet network distance with coordinates-based approaches,'' in {\em Proceedings. Twenty-First Annual Joint Conference of the IEEE Computer and Communications Societies}, vol.~1, pp.~170--179, IEEE, 2002.

\bibitem{gower1975generalized}
J.~C. Gower, ``Generalized procrustes analysis,'' {\em Psychometrika}, vol.~40, no.~1, pp.~33--51, 1975.

\bibitem{gower2004procrustes}
J.~C. Gower and G.~B. Dijksterhuis, {\em Procrustes problems}, vol.~30.
\newblock Oxford University Press, USA, 2004.

\bibitem{de2004sparse}
V.~De~Silva and J.~B. Tenenbaum, ``Sparse multidimensional scaling using landmark points,'' tech. rep., Stanford University, 2004.

\bibitem{platt2005fastmap}
J.~Platt, ``Fastmap, metricmap, and landmark mds are all nystr{\"o}m algorithms,'' in {\em International Workshop on Artificial Intelligence and Statistics}, pp.~261--268, PMLR, 2005.

\bibitem{delicado2024multidimensional}
P.~Delicado and C.~Pach{\'o}n-Garc{\'\i}a, ``Multidimensional scaling for big data,'' {\em Advances in Data Analysis and Classification}, pp.~1--22, 2024.

\bibitem{gower1968adding}
J.~C. Gower, ``Adding a point to vector diagrams in multivariate analysis,'' {\em Biometrika}, vol.~55, no.~3, pp.~582--585, 1968.

\bibitem{beck2008exact}
A.~Beck, P.~Stoica, and J.~Li, ``Exact and approximate solutions of source localization problems,'' {\em IEEE Transactions on signal processing}, vol.~56, no.~5, pp.~1770--1778, 2008.

\bibitem{borg2005modern}
I.~Borg and P.~J. Groenen, {\em Modern multidimensional scaling: Theory and applications}.
\newblock Springer Science \& Business Media, 2005.

\bibitem{so2007theory}
A.~M.-C. So and Y.~Ye, ``Theory of semidefinite programming for sensor network localization,'' {\em Mathematical Programming}, vol.~109, no.~2-3, pp.~367--384, 2007.

\bibitem{shi2023facial}
H.~Shi and Q.~Li, ``A facial reduction approach for the single source localization problem,'' {\em Journal of Global Optimization}, vol.~87, no.~2, pp.~831--855, 2023.

\bibitem{cox2000multidimensional}
T.~F. Cox and M.~A. Cox, {\em Multidimensional scaling}.
\newblock CRC press, 2000.

\bibitem{tenenbaum2000global}
J.~B. Tenenbaum, V.~d. Silva, and J.~C. Langford, ``A global geometric framework for nonlinear dimensionality reduction,'' {\em science}, vol.~290, no.~5500, pp.~2319--2323, 2000.

\bibitem{qi2013lagrangian}
H.-D. Qi, N.~Xiu, and X.~Yuan, ``A lagrangian dual approach to the single-source localization problem,'' {\em IEEE Transactions on Signal Processing}, vol.~61, no.~15, pp.~3815--3826, 2013.

\bibitem{kong2019classical}
L.~Kong, C.~Qi, and H.-D. Qi, ``Classical multidimensional scaling: A subspace perspective, over-denoising, and outlier detection,'' {\em IEEE Transactions on Signal Processing}, vol.~67, no.~14, pp.~3842--3857, 2019.

\bibitem{trosset2008out}
M.~W. Trosset and C.~E. Priebe, ``The out-of-sample problem for classical multidimensional scaling,'' {\em Computational statistics \& data analysis}, vol.~52, no.~10, pp.~4635--4642, 2008.

\bibitem{anderson2003generalized}
M.~J. Anderson and J.~Robinson, ``Generalized discriminant analysis based on distances,'' {\em Australian \& New Zealand Journal of Statistics}, vol.~45, no.~3, pp.~301--318, 2003.

\bibitem{gower1966some}
J.~C. Gower, ``Some distance properties of latent root and vector methods used in multivariate analysis,'' {\em Biometrika}, vol.~53, no.~3-4, pp.~325--338, 1966.

\bibitem{gower1982euclidean}
J.~C. Gower, ``Euclidean distance geometry,'' {\em Math. Sci.}, vol.~1, pp.~1--14, 1982.

\bibitem{more1993generalizations}
J.~J. Mor{\'e}, ``Generalizations of the trust region problem,'' {\em Optimization methods and Software}, vol.~2, no.~3-4, pp.~189--209, 1993.

\bibitem{larsson2025single}
M.~Larsson, V.~Larsson, K.~{\AA}str{\"o}m, and M.~Oskarsson, ``Single-source localization as an eigenvalue problem,'' {\em IEEE Transactions on Signal Processing}, 2025.

\bibitem{zhou2020robust}
S.~Zhou, N.~Xiu, and H.-D. Qi, ``Robust euclidean embedding via edm optimization,'' {\em Mathematical Programming Computation}, vol.~12, pp.~337--387, 2020.

\bibitem{drusvyatskiy2017noisy}
D.~Drusvyatskiy, N.~Krislock, Y.-L. Voronin, and H.~Wolkowicz, ``Noisy euclidean distance realization: robust facial reduction and the pareto frontier,'' {\em SIAM Journal on Optimization}, vol.~27, no.~4, pp.~2301--2331, 2017.

\bibitem{lichtenberg2024localization}
S.~Lichtenberg and A.~Tasissa, ``Localization from structured distance matrices via low-rank matrix recovery,'' {\em IEEE Transactions on Information Theory}, 2024.

\bibitem{Neidigh2002TrpCage}
J.~W. Neidigh, R.~M. Fesinmeyer, and N.~H. Andersen, ``Designing a 20‐residue protein,'' {\em Nature Structural Biology}, vol.~9, no.~7, pp.~425--430, 2002.

\bibitem{PDB1L2Y}
J.~W. Neidigh, R.~M. Fesinmeyer, and N.~H. Andersen, ``{1L2Y: NMR Structure of Trp‐Cage Miniprotein Construct TC5b}.'' Protein Data Bank (RCSB PDB), 2002.

\bibitem{pymol}
L.~Schrödinger and W.~DeLano, ``The pymol molecular graphics system.'' \url{http://www.pymol.org/pymol}, 2025.
\newblock Version 3.1.

\end{thebibliography}

\end{document}